\documentclass[acmtoms]{acmtrans2m}
\usepackage[rgb]{xcolor}
\usepackage[subrefformat=parens,labelformat=parens]{subfig}
\usepackage{caption}
\usepackage{moreverb, amsmath, amsfonts}
\usepackage{graphicx}
\usepackage{placeins}
\usepackage[ruled]{algorithm2e}
\usepackage{booktabs}
\usepackage{multirow}
\usepackage{xspace}
\usepackage{todonotes}

\usepackage{adjustbox}
\usepackage{scalefnt}
\usepackage{tikz}
\usetikzlibrary{patterns,calc}
\usepackage{pgf}

\usetikzlibrary{positioning}

\makeatletter

\newcommand{\lucbnc}[0]{SIESTA\_C\_BN\_1x1\xspace}
\newcommand{\ludna}[0]{DNA\_16\xspace}
\newcommand{\dgacpn}[0]{DGDFT\_ACPNR4\_60\xspace}

\newenvironment{topbox}[2]{%
    \pgfmathsetlength\@tempdima{#2}%
    \pgfmathsetlength\pgf@yc{\pgfkeysvalueof{/pgf/inner ysep}}%
    \advance\@tempdima by -2\pgf@yc
    \begin{lrbox}{\@tempboxa}%
        \begin{minipage}[t]{#1}%
           \vspace{0pt}%
}{%
        \end{minipage}%
    \end{lrbox}%
    \ifdim\@tempdima>\dp\@tempboxa
        \dp\@tempboxa=\@tempdima
    \fi
    \box\@tempboxa
}

\makeatother

\tikzstyle{block} = [rectangle,fill=none,solid, draw=black, text centered,minimum width=2em, minimum height=2em,inner sep=0pt]
\DeclareRobustCommand\circled[1]{\tikz[baseline=(char.base)]{ \node[shape=circle, fill=white,draw,inner sep=0pt,text width= 1 em ,minimum size = 1 em, text centered] (char) {#1};}}
\tikzstyle{treenode} = [rectangle,fill=white,solid, draw=black, text centered,minimum size=1.5em, text width=1.5em,inner sep=0pt]
\DeclareRobustCommand\supernode[1]{\tikz[baseline=(char.base)]{ \node[shape=rectangle, fill=white,draw,inner sep=0pt,text width= 1 em ,minimum size = 1 em, text centered] (char) {#1};}}

\newcommand{\elm}[2]{\fill[fill=gray!40] ($(#1 em,-#2 em )+(0.5 em,-0.5 em)$) circle [radius=0.4 em]}
\newcommand{\welm}[2]{\fill[fill=white] ($(#1 em,-#2 em )+(0.5 em,-0.5 em)$) circle [radius=0.45 em]}
\newcommand{\felm}[2]{\fill[fill=gray!80] ($(#1 em,-#2 em )+(0.15 em,-0.15 em)$) rectangle +(0.7em,-0.7em)}
\newcommand{\wfelm}[2]{\fill[fill=white] ($(#1 em,-#2 em )+(0.1 em,-0.1 em)$) rectangle +(0.8em,-0.8em)}

\newcommand{\edblk}[3]{\node[block,minimum width=#2 em, minimum height=#2 em] at (#3 em,-#3 em) {#1};}

\newcommand{\dblk}[3]{
    \pgfmathtruncatemacro{\start}{#3 -0.5*#2};
    \pgfmathtruncatemacro{\stop}{#3 +0.5*#2-1};
  \foreach \i in {\start,...,\stop}{
    \foreach \j in {\i,...,\stop}{
      \elm{\i}{\j};
    }
  }
\edblk{\supernode{#1}}{#2}{#3};
}

\newcommand{\mapblock}[4]{
    \pgfmathtruncatemacro{\curproc}{ int(mod(#2-1,#3))*#4 + int(mod(#1-1,#4)) +1 };
}

\newcommand{\colorproc}[3]{
  \pgfmathsetmacro{\h}{#1 /(#2*#3)}
  \definecolor{proccolor}{hsb}{\h, 0.5, 0.8}
}

\newcommand{\blk}[6]{
\node[block,minimum width=#2em, minimum height=#3em #6] at (#4em,-#5em) {#1};
}

\newcommand{\dblknew}[4]{
\blk{#1}{#2}{#2}{#3}{#3}{#4}
}

\newcommand{\mdblk}[6]{
    \pgfmathtruncatemacro{\x}{ int(mod(#2-1,#5))*#6 + int(mod(#2-1,#6)) +1 };
 \pgfmathsetmacro{\h}{\x /(#5*#6)}
\definecolor{currentcolor}{hsb}{\h, 0.5, 0.8}
\node[block,minimum width=#3 em, minimum height=#3 em,fill=currentcolor] at (#4 em,-#4 em) {\Proc{\x}};
}

\newcommand{\elublk}[5]{
\node[block,minimum width=#2em, minimum height=#3em] at (#4em,-#5em) {#1};\node[block,minimum width=#3em, minimum height=#2em] at (#5em,-#4em) {#1};
}

\newcommand{\elublkFill}[5]{
\node[block,draw=none,minimum width=#2em, minimum height=#3em,fill=white,opacity=1] at (#4em,-#5em) {#1};
\draw[draw=black] ($(#4em,-#5em)- 0.5*(#2em,-#3em)$) -- ($(#4em,-#5em)- 0.5*(#2em,#3em)$);
\draw[draw=black] ($(#4em,-#5em) + 0.5*(#2em,#3em)$) -- ($(#4em,-#5em) + 0.5*(#2em,-#3em)$);

\node[block,draw=none,minimum width=#3em, minimum height=#2em,fill=white,opacity=1] at (#5em,-#4em) {#1};
\draw[draw=black] ($(#5em,-#4em)- 0.5*(#3em,-#2em)$) -- ($(#5em,-#4em)- 0.5*(-#3em,-#2em)$);
\draw[draw=black] ($(#5em,-#4em)- 0.5*(#3em,#2em)$) -- ($(#5em,-#4em)- 0.5*(-#3em,#2em)$);

}

\newcommand{\mlublk}[9]{
    \pgfmathtruncatemacro{\x}{ int(mod(#3-1,#8))*#9 + int(mod(#2-1,#9)) +1 };
 \pgfmathsetmacro{\h}{\x /(#8*#9)}
\definecolor{currentcolor}{hsb}{\h, 0.5, 0.8}
 \node[block,minimum width=#4em, minimum height=#5em,fill=currentcolor] at (#6em,-#7em) {\Proc{\x}};

    \pgfmathtruncatemacro{\x}{ int(mod(#2-1,#8))*#9 + int(mod(#3-1,#9)) +1 };
 \pgfmathsetmacro{\h}{\x /(#8*#9)}
\definecolor{currentcolor}{hsb}{\h, 0.5, 0.8}
 \node[block,minimum width=#5em, minimum height=#4em,fill=currentcolor] at (#7em,-#6em) {\Proc{\x}};
}

\newcommand{\lublk}[5]{
    \pgfmathtruncatemacro{\starti}{#4 -0.5*#2};
    \pgfmathtruncatemacro{\stopi}{#4 +0.5*#2-1};
    \pgfmathtruncatemacro{\startj}{#5 -0.5*#3};
    \pgfmathtruncatemacro{\stopj}{#5 +0.5*#3-1};
  \foreach \i in {\starti,...,\stopi}{
    \foreach \j in {\startj,...,\stopj}{
      \elm{\i}{\j};
    }
  }
\elublk{#1}{#2}{#3}{#4}{#5};
}

\newcommand{\flublk}[5]{
    \pgfmathtruncatemacro{\starti}{#4 -0.5*#2};
    \pgfmathtruncatemacro{\stopi}{#4 +0.5*#2-1};
    \pgfmathtruncatemacro{\startj}{#5 -0.5*#3};
    \pgfmathtruncatemacro{\stopj}{#5 +0.5*#3-1};
  \foreach \i in {\starti,...,\stopi}{
    \foreach \j in {\startj,...,\stopj}{
      \felm{\i}{\j};
    }
  }
\elublk{#1}{#2}{#3}{#4}{#5};
}

\newcommand{\ignore}[1]{}

\newcommand{\Or}{\ensuremath{\mathcal{O}}\xspace}

\newcommand{\Tr}{\ensuremath{\mathrm{Tr}}\xspace}

\newcommand{\CS}{\ensuremath{\mathcal C}\xspace}
\newcommand{\JS}{\ensuremath{\mathcal J}\xspace}
\newcommand{\IS}{\ensuremath{\mathcal I}\xspace}
\newcommand{\KS}{\ensuremath{\mathcal K}\xspace}
\newcommand{\RS}{\ensuremath{\mathcal R}\xspace}

\newcommand{\Amat}{\ensuremath{A}\xspace}
\newcommand{\Lmat}{\ensuremath{L}\xspace}

\newcommand{\Proc}[1]{\ensuremath{P_{#1}}\xspace}

\newcommand{\etc}{\textit{etc}.\xspace}

\newcommand{\ie}{\textit{i.e.}\xspace}

\newcommand{\abs}[1]{\lvert#1\rvert} 
\newcommand{\norm}[1]{\lVert#1\rVert}

\newcommand{\LL}[1]{{\color{red}~\textsf{[LL: #1]}}}

\newcommand{\REV}[1]{{#1}}
\newcommand{\REVTWO}[1]{{#1}}

\newcommand{\pselinv}{\texttt{PSelInv}\xspace}

\newcommand{\pardiso}{\texttt{PARDISO}\xspace}
\newcommand{\mumps}{\texttt{MUMPS}\xspace}
\newcommand{\superlu}{\texttt{SuperLU\_DIST}\xspace}
\newcommand{\parmetis}{\texttt{ParMETIS}\xspace}
\newcommand{\ptscotch}{\texttt{PT-Scotch}\xspace}

\title{\texttt{PSelInv} -- A Distributed Memory Parallel Algorithm for
Selected Inversion : the Symmetric Case} 
\author{
Mathias Jacquelin, Lawrence Berkeley National Laboratory\\ 
Lin Lin, University of California, Berkeley and Lawrence Berkeley National
Laboratory \\
Chao Yang, Lawrence Berkeley National Laboratory
}

\begin{abstract}
We describe an efficient parallel implementation of the selected
inversion algorithm for distributed memory computer systems, which we
call \texttt{PSelInv}.  The \texttt{PSelInv} method computes selected
elements of a general sparse matrix $A$ that can be decomposed as $A =
LU$, where $L$ is lower triangular and $U$ is upper triangular.  The
implementation described in this paper focuses on the case of sparse
symmetric matrices.  It contains an interface that is compatible with
the distributed memory parallel sparse direct factorization
\texttt{SuperLU\_DIST}.  However, the underlying data structure and
design of \texttt{PSelInv} allows it to be easily combined with other
factorization routines such as \texttt{PARDISO}.  We discuss general
parallelization strategies such as data and task distribution schemes.
In particular, we describe how to exploit the concurrency exposed by the
elimination tree associated with the $LU$ factorization of $A$.  We
demonstrate the efficiency and accuracy of \texttt{PSelInv} by
presenting a number of numerical experiments.  In particular, we show
that \texttt{PSelInv} can run efficiently on more than $4,000$
\REVTWO{cores}
for a modestly sized matrix.  We also demonstrate how \texttt{PSelInv}
can be used to accelerate large-scale electronic structure calculations.
\end{abstract}






\begin{document}

\maketitle

\nocite{LinLuYingE2009,LinLuYingCarE2009,LinYangLuEtAl2011,LinYangMezaEtAl2011}

\section{Introduction}\label{sec:intro}

Let $A\in \mathbb{C}^{N\times N}$ be a non-singular sparse matrix.  We
are interested in computing {\em selected elements} of $A^{-1}$, defined
as
\begin{equation}
  \{(A^{-1})_{i,j}\vert \ \ \mbox{for} \ \ 1\le i,j\le N, \ \ 
  \mbox{such that} \ \ A_{i,j}\ne 0 \}.
  \label{eqn:selelem}
\end{equation}
Sometimes, we only need to compute a subset of these selected elements,
for example, the diagonal elements of $A^{-1}$. \REV{The most
straightforward way to obtain these selected elements of $A^{-1}$ is to
compute the full inverse of $A$ and then extract the selected elements.
But this is often prohibitively expensive in practice.  
It turns out that in order to compute these selected elements of $A^{-1}$,
some additional elements of $A^{-1}$ often need to be computed. However,
the overall set of nonzero elements that need to be computed 
often remains a small percentage of all elements of $A^{-1}$ due
to the sparsity structure of $A$.}

The selected elements of $A^{-1}$ defined by \eqref{eqn:selelem} can be 
used to obtain trace estimation of the form
\begin{equation}
  \Tr[A^{-1}]\quad \mbox{or} \quad \Tr[A^{-1}B^{T}],
  \label{eqn:traceest}
\end{equation}
if the sparsity pattern of $B\in\mathbb{C}^{N\times N}$ is contained in
the sparsity pattern of $A$, i.e. $\{(i,j)\vert B_{i,j}\ne 0\}\subset
\{(i,j)\vert A_{i,j}\ne 0\}$.  The computation of selected elements of $A^{-1}$,
together with the trace estimation of the form \eqref{eqn:traceest} 
arise in a number of scientific
computing applications including density functional
theory (DFT)~\cite{HohenbergKohn1964,KohnSham1965}, dynamical mean field
theory (DMFT)~\cite{KotliarSavrasovHauleEtAl2006},
Poisson-Boltzmann equation~\cite{XuMaggs2013}, and uncertainty
quantification~\cite{BekasCurioniFedulova2009} etc.. 

It is possible to compute selected elements of $A^{-1}$ 
by iterative methods such as the Lanczos 
algorithm~\cite{Lanczos1950,SidjeSaad2011}, combined with 
Monte Carlo~\cite{BekasKokiopoulouSaad2007} or deterministic probing
techniques~\cite{TangSaad2012}.  This type of methods work well
if $A^{-1}$ is a banded matrix, or becomes a banded matrix after 
elements with absolute value less than $\varepsilon$ have been truncated,
and a sparse factorization of $A$ is prohibitively expensive
to perform.

In this paper we focus on using a sparse direct method to
compute selected elements of $A^{-1}$. We assume that a sparse 
LU factorization (or $LDL^{T}$ factorization if $A$ is symmetric) 
of $A$ is computationally feasible. The main advantage of a direct
method is that we do not need to make assumptions on
the decay property of $A^{-1}$.  \REV{The disadvantage of this direct
method is that it actually computes a \textit{superset} of the selected
elements of $A^{-1}$ as defined in Eq.~\eqref{eqn:selelem}.  In
particular, all elements of $A^{-1}$ indexed by the union of the
nonzero index sets of the $L$ and $U$ factors need to be computed.  
It should be noted that as long as $L$ and $U$ remain sparse, 
computing elements of $A^{-1}$ restricted to this superset can still be
much faster than computing the full inverse.}

Sparse direct methods for computing selected elements of $A^{-1}$ were
first proposed in the
papers~\cite{TakahashiFaganChin1973,ErismanTinney1975}.  The use of
elimination tree for \REV{computing selected elements of inverse} was presented
in~\cite{CampbellDavis1995} \REV{for a sequential algorithm}.  
Motivated by quantum transport simulations, Li et al.
\cite{LiAhmedKlimeckDarve2008,LiDarve2012,LiWuDarve2013} developed the Fast Inverse
using Nested Dissection (FIND) algorithm for computing the diagonal of
$A^{-1}$. \REV{Some related work has recently been described in
~\cite{CauleyBalakrishnanKlimeckEtAl2012,EastwoodWan2013}.
The FIND algorithm is in principle applicable to matrices
with a general sparsity pattern, but so far its implementation focuses on
structured matrices obtained from second-order partial differential
operators discretized by a finite difference scheme. The
implementation of FIND is not yet publicly available. Nor is it scalable to a
large number of processors. }
\REV{[Lin et al. 2009b]} developed the Hierarchically Schur Complement
(HSC) method for a similar type of matrix arising from density
functional theory based electronic structure calculations.
The method has also been generalized and applied to quantum transport
calculations~\cite{HetmaniukZhaoAnantram2013}. 

An efficient implementation of the selected inversion algorithm
for a general symmetric matrix, called \texttt{SelInv}, was presented 
in~\REV{[Lin et al. 2011b]}, and is publicly available. 
Amestoy et
al.~\REVTWO{\cite{AmestoyDuffLExcellentEtAl2012a}}
considered a more general parallel matrix inversion method for computing any
subset of entries of $A^{-1}$.  They implemented their algorithm in the
\mumps package~\cite{mumps}, which is based on the multifrontal method.
In this algorithm, the actual set of computed entries of $A^{-1}$ 
\REVTWO{contains entries on the critical path of the requested
entries to the root of the elimination tree, and therefore this}
is also a \textit{superset} of the requested entries of $A^{-1}$.  This
method is more efficient than our algorithm when a small number of entries of
$A^{-1}$ are requested.  However, when a relatively large number of
entries are to be computed such as in the computation of the selected
elements defined by~\eqref{eqn:selelem}, our algorithm can 
reuse more efficiently the information shared among different entries of
$A^{-1}$, and our numerical results indicate that our algorithm is more
efficient than the parallel matrix inversion method implemented in
\mumps.
In addition to the work presented
in~\cite{AmestoyDuffLExcellentEtAl2012,AmestoyDuffLExcellentEtAl2012a},
parallel implementation of algorithms for computing selected elements of
inverse tailored to
matrices obtained from a finite difference discretization of a
second-order partial differential operator have appeared in a number of
publications~\cite{PetersenLiStokbroEtAl2009,LinYangLuEtAl2011}.  The
method developed by Petersen et al.~\cite{PetersenLiStokbroEtAl2009},
which was designed for quasi-1D quantum transport problem, is scalable
to a relatively small ($32\sim 64$) number of processors.
In~\REV{[Lin et al. 2011a]}, we described an efficient parallel
implementation for discretized 2D Laplacian type of operators, and
demonstrated the efficiency of the implementation when it was used to
solve a problem with billions degrees of freedom on 4,096 processors.
However, this implementation cannot be used to perform a selected
inversion of a general symmetric sparse matrix with an arbitrary
sparsity pattern. 


The purpose of this paper is to extend the selected inversion algorithm presented 
in~\REV{[Lin et al. 2011b]} and describe an implementation of a
parallel selected inversion algorithm designed for distributed 
memory parallel computers. Such an implementation allows us
to solve much larger problems by utilizing more computational 
resources available on high performance computers. \REV{The present work
is more general than the previous work in [Lin et al. 2011a] which
assumes a balanced binary elimination tree and is only applicable to
structured sparse matrices obtained from finite difference discretization 
of differential operators. It can utilize far more number of processors
than the sequential algorithm described in [Lin et al. 2011b] for
general sparse matrices.}

The parallel implementation of the selected inversion algorithm 
we present in this paper uses a more general data distribution
and communication strategy to divide the work among a large 
number of processors to achieve multiple levels of concurrency.
We name our implementation \pselinv and the software is publicly
available. 
It is publicly available\footnote{\url{http://www.pexsi.org/}, distributed
under the BSD license}. 
Our first implementation focuses on the case of sparse symmetric matrices.  In principle, the \pselinv package can be 
interfaced with any sparse $LU$ and $LDL^{T}$ factorization routines. 
Our current implementation provides an interface to the
\superlu~\cite{LiDemmel2003} package.  The user has the option 
of using either the \parmetis~\cite{KarypisKumar1998a} software or the
\ptscotch~\cite{ChevalierPellegrini2008} package to reorder the matrix
in parallel to minimize the non-zero fill in the sparse matrix factors.

The rest of the paper is organized as follows.  We review the basic
idea of the selected inversion method in
Section~\ref{sec:selinv}, and discuss various implementation issues for
the distributed memory parallel selected inversion algorithm in
Section~\ref{sec:parallelization}.  The numerical results with
applications to various matrices from including Harwell-Boeing Test Collection
\cite{HarwellBoeing}, the University of Florida Matrix
Collection\cite{FloridaMatrix}, and also applications from density
functional theory are given in Section~\ref{sec:numerical}, followed by
the conclusion and the future work discussion in Section~\ref{sec:conclusion}.

Standard linear algebra notation is used for vectors and matrices
throughout the paper.  We use $A_{i,j}$ to denote the $(i,j)$-th entry
of the matrix $A$, and $f_{i}$ to denote the $i$-th entry of the vector
$f$. With slight abuse of notation, both a supernode index and the set
of column indices associated with a supernode are denoted by uppercase
script letters such as $\IS,\JS,\KS$ \etc.  Furthermore, we use
$A_{i,*}$ and $A_{*,j}$ to denote the $i$-th row and the $j$-th column
of $A$, respectively. Similarly, $A_{\IS,*}$ and $A_{*,\JS}$ are used to
denote the $\IS$-th block row and the $\JS$-th block column of $A$,
respectively. \REV{$A_{\IS,\JS}^{-1}$ denotes the $(\IS,\JS)$-th block of the
matrix $A^{-1}$, i.e. $A_{\IS,\JS}^{-1}\equiv (A^{-1})_{\IS,\JS}$.  When
the block $A_{\IS,\JS}$ itself is invertible, its inverse is denoted by
$(A_{\IS,\JS})^{-1}$ to distinguish from $A_{\IS,\JS}^{-1}$.}

\section{Selected inversion algorithm}\label{sec:selinv}

\subsection{Basic formulation}\label{subsec:basic}

Although this paper focuses on the computation of selected elements of
$A^{-1}$ when $A$ is a sparse symmetric matrix, the idea of the selected
inversion algorithm can be given for a general square matrix
$A$, \REV{and will be used in our ongoing work for computing the
selected elements of $A^{-1}$ for asymmetric matrices.}
The standard approach for computing $A^{-1}$ is to first decompose
the general matrix $A$ using the LU factorization 
\begin{equation}
	A = LU
	\label{eqn:LU}
\end{equation}
where $L$ is a unit lower triangular matrix and $U$ is an upper
triangular matrix.  In order to stabilize the computation, matrix
reordering and partial pivoting~\cite{GolubVan1996} are usually applied to the
matrix of A, and the general form of the LU factorization can be
given as
\begin{equation}
	PAQ = LU
	\label{eqn:permuteLU}
\end{equation}
where $P$ and $Q$ are two permutation matrices.  To simplify the
discussion below we use Eq.~\eqref{eqn:LU} and assume $A$ has already been
permuted.

Given the $LU$
factorization, the most straightforward way to compute selected elements
of $A^{-1}$ is to obtain $A^{-1}\equiv(x_{1},x_{2},\ldots,x_{n})$ by solving a
number of triangular systems
\begin{equation}
	L y_{j}=e_{j}, \quad U x_{j}=y_{j}.
	\label{eqn:triangular}
\end{equation}
for $j=1,2,\ldots,n$, and $e_{j}$ is the $j$-th column of the $n\times
n$ identity matrix. Such a procedure, which will be referred to as the
direct inversion algorithm, is generally very costly even when
$A$ is sparse. The direct inversion algorithm performs too much computation
when only a small number of selected elements of the
inverse matrix are needed.  

An alternative algorithm is the selected inversion algorithm, which
accurately computes all the selected elements of $A^{-1}$.  The idea of
the selected inversion method  originates
from~\cite{TakahashiFaganChin1973,ErismanTinney1975}, and the algorithm
and its variants have been discussed in a number of
recent
works~\REV{[Lin et al. 2009b; 2011b],
\cite{AmestoyDuffLExcellentEtAl2012,LiAhmedKlimeckDarve2008,LiDarve2012,KuzminLuisierSchenk2013}}.  
The selected inversion algorithm can be understood as follows. We
first partition the matrix $A$ into $2\times 2$
blocks of the form
\begin{equation}
	A = \begin{pmatrix}
		A_{1,1} & A_{1,2}\\
		A_{2,1} & A_{2,2}
	\end{pmatrix},
	\label{}
\end{equation}
where $A_{1,1}$ is a scalar of size $1\times 1$.  We can
write $A_{1,1}$ as a product of two scalars $L_{1,1}$ and $U_{1,1}$.
In particular, we can pick $L_{1,1}=1$ and $U_{1,1}=A_{1,1}$. Then 
\begin{equation}
	A = \begin{pmatrix}
		L_{1,1} & 0\\
		L_{2,1} & I 
	\end{pmatrix}
	\begin{pmatrix}
		U_{1,1} & U_{1,2}\\
		0      & S_{2,2}
	\end{pmatrix}
	\label{eqn:LU2by2}
\end{equation}
where
\REV{
\begin{equation}
	L_{2,1}=A_{2,1} (U_{1,1})^{-1}, \quad U_{1,2} = (L_{1,1})^{-1}	A_{1,2}.
	\label{}
\end{equation}
}
The $L$ and $U$ factors are usually directly accessible in a standard
$LU$ factorization, and
\begin{equation}
	S_{2,2} = A_{2,2} - L_{2,1} U_{1,2}	
	\label{}
\end{equation}
is the Schur complement.  Using the decomposition given by 
Eq.~\eqref{eqn:LU2by2}, we can express $A^{-1}$ as
\REV{
\begin{equation}
	A^{-1} = \begin{pmatrix}
		(U_{1,1})^{-1} (L_{1,1})^{-1} + (U_{1,1})^{-1} U_{1,2} S^{-1}_{2,2}
		L_{2,1} (L_{1,1})^{-1} & - (U_{1,1})^{-1} U_{1,2} S^{-1}_{2,2} \\
		-S^{-1}_{2,2} L_{2,1} (L_{1,1})^{-1} & S_{2,2}^{-1}
	\end{pmatrix}.
	\label{eqn:Ainv2by2}
\end{equation}
}
\REV{Since $S_{2,2}$ is the same as $S$ here, without
ambiguity $S_{2,2}^{-1}\equiv (S^{-1})_{2,2}$ can be used.}
To simplify the notation, we define the normalized $LU$ factors as
\REV{
\begin{equation}
	\hat{L}_{1,1} = L_{1,1}, \quad \hat{U}_{1,1}=U_{1,1},\quad
	\hat{L}_{2,1} = L_{2,1}(L_{1,1})^{-1}, \quad \hat{U}_{1,2} = (U_{1,1})^{-1}
	U_{1,2},
	\label{}
\end{equation}
}
and Eq.~\eqref{eqn:Ainv2by2} can be equivalently given by 
\REV{
\begin{equation}
	A^{-1} = \begin{pmatrix}
    (\hat{U}_{1,1})^{-1} (\hat{L}_{1,1})^{-1} + \hat{U}_{1,2}
    S_{2,2}^{-1}
		\hat{L}_{2,1} & - \hat{U}_{1,2} S_{2,2}^{-1} \\
		-S_{2,2}^{-1} \hat{L}_{2,1}  & S_{2,2}^{-1}
	\end{pmatrix}.
	\label{eqn:Ainv2by2normal}
\end{equation}
}
Let us denote by $\CS$ the set of indices
\begin{equation}
	\{i|\left(L_{2,1}\right)_{i} \ne 0\} \cup 
	\{j|\left(U_{1,2}\right)_{j} \ne 0\},
	\label{}
\end{equation} 
and assume $S_{2,2}^{-1}$ has already been computed. 
From Eq.~\eqref{eqn:Ainv2by2normal} it can be readily observed that
in order to compute the \REVTWO{$i$-th element of
$A^{-1}_{2,1} \equiv -S^{-1}_{2,2} \hat{L}_{2,1}$} for $i\in \CS$,  
we only need the entries 
\begin{equation}
\left\{\left( S_{2,2}^{-1} \right)_{i,j} | i\in \CS, j\in \CS\right\}.
	\label{eqn:selectentry2x2}
\end{equation}
The same set of entries of $S_{2,2}^{-1}$ are required to compute 
selected entries of \REV{$A_{1,2}^{-1} \equiv -\hat{U}_{1,2}
S^{-1}_{2,2}$}. 
No additional entries of $S_{2,2}^{-1}$ are needed to complete
the computation of \REV{$A_{1,1}^{-1}$}, which involves the matrix product 
of selected entries of $\hat{U}_{1,2}$ and \REV{$A_{2,1}^{-1}$}.
This procedure can be repeated recursively to compute selected elements
of $S_{2,2}^{-1}$ until $S_{2,2}$ is a scalar of size $1$. A pseudo-code for
demonstrating this column-based selected inversion algorithm for symmetric matrix is
given in~\cite{LinYangMezaEtAl2011}.

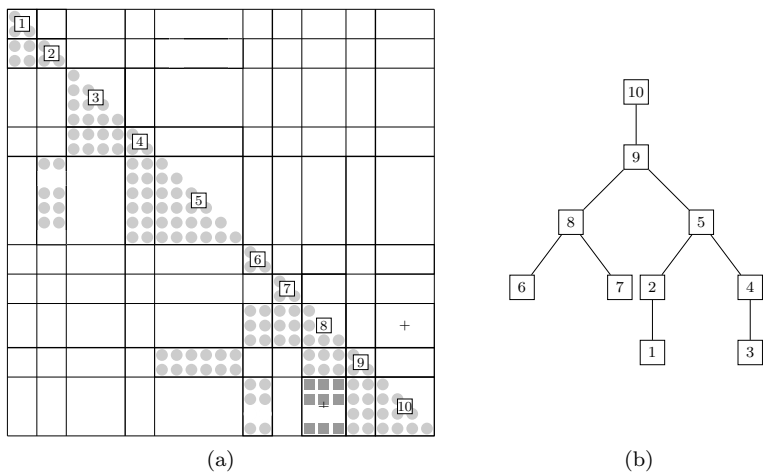
\begin{figure}
\centering
\subfloat[]{\label{fig:matrixA.matrix}
\begin{adjustbox}{width=.45\textwidth}
\begin{tikzpicture}

\dblk{1}{2}{1};
\dblk{2}{2}{3};
\dblk{3}{4}{6};
\dblk{4}{2}{9};
\dblk{5}{6}{13};
\dblk{6}{2}{17};
\dblk{7}{2}{19};
\dblk{8}{3}{21.5};
\dblk{9}{2}{24};
\dblk{10}{4}{27};

\draw (0em,0em) --   (29em ,0em);
\draw (0em,-2em) --  (29em ,-2em);
\draw (0em,-4em) --  (29em ,-4em);
\draw (0em,-8em) --  (29em ,-8em);
\draw[thick] (0em,-10em) -- (29em ,-10em);
\draw (0em,-16em) -- (29em ,-16em);
\draw (0em,-18em) -- (29em ,-18em);
\draw (0em,-20em) -- (29em ,-20em);
\draw[thick] (0em,-23em) -- (29em ,-23em);
\draw (0em,-25em) -- (29em ,-25em);
\draw (0em,-29em) -- (29em ,-29em);

\draw (0em,0em) --  (0em ,-29em);
\draw (2em ,0em) -- (2em ,-29em);
\draw (4em ,0em) -- (4em ,-29em);
\draw[thick] (8em ,0em) -- (8em ,-29em);
\draw (10em,0em) -- (10em,-29em);
\draw (16em,0em) -- (16em,-29em);
\draw[thick] (18em,0em) -- (18em,-29em);
\draw (20em,0em) -- (20em,-29em);
\draw (23em,0em) -- (23em,-29em);
\draw[thick] (25em,0em) -- (25em,-29em);
\draw (29em,0em) -- (29em,-29em);

\lublk{}{2}{2}{1}{3};
\lublk{}{2}{6}{3}{13};

\elublkFill{}{2}{1}{3}{11.5};
\elublkFill{}{2}{1}{3}{15.5};
\draw[draw=black] ($(15.5em,-3em) + 0.5*(1em,2em)$) -- ($(15.5em,-3em) + 0.5*(1em,-2em)$);
\draw[draw=black] ($(3em,-15.5em) + 0.5*(2em,-1em)$) -- ($(3em,-15.5em) + 0.5*(-2em,-1em)$);

\lublk{}{4}{2}{6}{9};
\lublk{}{2}{6}{9}{13};


%
\lublk{}{6}{2}{13}{24};
%

\lublk{}{2}{3}{17}{21.5};
\lublk{}{2}{4}{17}{27};


      \welm{16}{27};
      \welm{17}{27};

\lublk{}{2}{3}{19}{21.5};

\lublk{}{3}{2}{21.5}{24};
\flublk{+}{3}{4}{21.5}{27};

      \wfelm{20}{27};
      \wfelm{21}{27};
      \wfelm{22}{27};
\lublk{}{2}{4}{24}{27};

%
%
%

\end{tikzpicture}
\end{adjustbox}
}
\quad
\subfloat[]{\label{fig:matrixA.etree}
\begin{adjustbox}{width=.35\textwidth}
\begin{tikzpicture}

\node[treenode] (10) at (0,-4em) {10};
\node[treenode] (9) at (0,-8em) {9};

\node[treenode] (8) at (-4 em,-12em) {8};
\node[treenode] (5) at (4 em,-12em) {5};

\node[treenode] (6) at (-7 em,-16em) {6};
\node[treenode] (7) at (-1 em,-16em) {7};
\node[treenode] (2) at (1 em,-16em) {2};
\node[treenode] (4) at (7 em,-16em) {4};

\node[treenode] (1) at (1 em,-20em) {1};
\node[treenode] (3) at (7 em,-20em) {3};

\draw (10) -- (9);
\draw (9) -- (5);
\draw (9) -- (8);
\draw (5) -- (2);
\draw (5) -- (4);
\draw (8) -- (6);
\draw (8) -- (7);

\draw (2) -- (1);
\draw (4) -- (3);

\begin{scope}[xscale = 1.2, yscale=1.3, shift={(-2.5,-0.1)}]
\begin{scope}[yscale=-1, x=0.5 em, y=0.5 em, shift={(-533.625,-648.65625)}]
  \path[draw=none,line join=miter,line cap=butt,line width=0.081pt]
    (539.2699,653.8548) .. controls (539.8629,653.3332) and (540.4902,652.8506) ..
    (541.1464,652.4113) .. controls (542.5908,651.4445) and (544.2441,650.6722) ..
    (545.9821,650.6792) .. controls (547.1051,650.6842) and (548.1973,651.0120) ..
    (549.2955,651.2468) .. controls (550.3936,651.4817) and (551.5598,651.6205) ..
    (552.6222,651.2568) .. controls (553.7811,650.8598) and (554.7090,649.8476) ..
    (555.0039,648.6586);
  \path[draw=none,line join=miter,line cap=butt,line width=0.081pt]
    (539.9195,661.0723) .. controls (540.5169,659.7651) and (541.5296,658.6516) ..
    (542.7743,657.9330) .. controls (544.0190,657.2144) and (545.4897,656.8943) ..
    (546.9204,657.0305) .. controls (548.5543,657.1861) and (550.1149,657.9154) ..
    (551.7561,657.8966) .. controls (552.9762,657.8826) and (554.1869,657.4364) ..
    (555.1243,656.6552) .. controls (556.0616,655.8741) and (556.7189,654.7637) ..
    (556.9526,653.5661);
  \path[draw=none,line join=miter,line cap=butt,line width=0.081pt]
    (535.8777,665.4749) .. controls (539.6178,663.8277) and (543.6998,662.9618) ..
    (547.7865,662.9488) .. controls (548.6527,662.9458) and (549.5192,662.9811) ..
    (550.3848,662.9488) .. controls (553.6313,662.8277) and (556.8355,661.7342) ..
    (559.4787,659.8453);
  \path[draw=none,line join=miter,line cap=butt,line width=0.081pt]
    (533.7125,671.7541) .. controls (534.5027,670.2971) and (535.9783,669.2313) ..
    (537.6099,668.9393) .. controls (538.8369,668.7196) and (540.0952,668.9183) ..
    (541.3284,669.0996) .. controls (542.5617,669.2810) and (543.8311,669.4438) ..
    (545.0439,669.1558) .. controls (546.5353,668.8016) and (547.7803,667.8073) ..
    (549.1578,667.1349) .. controls (550.7532,666.3562) and (552.5295,666.0100) ..
    (554.1378,665.2584) .. controls (554.8790,664.9120) and (555.5780,664.4825) ..
    (556.3031,664.1036) .. controls (557.7781,663.3327) and (559.3628,662.7720) ..
    (560.9944,662.4436);
  \path[draw=none,line join=miter,line cap=butt,line width=0.081pt]
    (533.6403,675.7237) .. controls (535.0555,675.8347) and (536.4942,675.6363) ..
    (537.8264,675.1463) .. controls (540.5119,674.1585) and (542.6494,672.0585) ..
    (545.2604,670.8880) .. controls (545.8442,670.6262) and (546.4491,670.4124) ..
    (547.0647,670.2384) .. controls (550.0786,669.3865) and (553.2715,669.5028) ..
    (556.3752,669.0836) .. controls (559.7657,668.6256) and (563.0670,667.5170) ..
    (566.0466,665.8358);
  \path[draw=none,line join=miter,line cap=butt,line width=0.081pt]
    (535.8056,682.4359) .. controls (542.9361,686.6844) and (542.0366,676.1286) ..
    (546.3430,675.8680) .. controls (547.3991,675.8041) and (548.4586,675.8690) ..
    (549.5163,675.8445) .. controls (550.5740,675.8200) and (551.6472,675.7014) ..
    (552.6222,675.2906) .. controls (553.5736,674.8897) and (554.3908,674.2281) ..
    (555.1507,673.5292) .. controls (555.9107,672.8304) and (556.6308,672.0833) ..
    (557.4579,671.4654) .. controls (558.7786,670.4786) and (560.3595,669.8434) ..
    (561.9958,669.6420) .. controls (563.6322,669.4406) and (565.3199,669.6735) ..
    (566.8405,670.3106);
  \path[draw=none,line join=miter,line cap=butt,line width=0.081pt]
    (542.4456,684.2402) .. controls (543.1264,682.5544) and (544.7014,681.2606) ..
    (546.4874,680.9202) .. controls (546.9627,680.8296) and (547.4476,680.8026) ..
    (547.9308,680.7759) .. controls (548.6531,680.7360) and (549.3756,680.6962) ..
    (550.0961,680.6315) .. controls (550.6254,680.5840) and (551.1536,680.5230) ..
    (551.6839,680.4872) .. controls (552.3847,680.4398) and (553.0942,680.4356) ..
    (553.7770,680.2706) .. controls (554.9290,679.9924) and (555.9250,679.2768) ..
    (556.8217,678.5018) .. controls (557.7184,677.7268) and (558.5554,676.8720) ..
    (559.5509,676.2289) .. controls (560.7505,675.4539) and (562.1618,675.0104) ..
    (563.5890,674.9599) .. controls (565.0162,674.9094) and (566.4553,675.2520) ..
    (567.7066,675.9402);
\end{scope}
\end{scope}

\end{tikzpicture}

\end{adjustbox}
\vspace{1cm}
}
\caption{\protect\subref{fig:matrixA.matrix} A structurally symmetric matrix $A$ of size $29\times 29$
divided into $10$ supernodes.  The nonzero matrix elements in $A$ are
labeled by round dots and the extra fill-in elements in $L$ are labeled
by squares.~\protect\subref{fig:matrixA.etree} \REV{The elimination tree corresponding to} the matrix
$A$ and its supernode partitioning.} 
\label{fig:matrixA}
\end{figure}

In practice, a column-based sparse factorization and selected
inversion algorithm may not be efficient due to the lack of level 3 BLAS 
operations.  For a sparse matrix $A$, 
the columns of $A$ and the $L$ factor can be
partitioned into supernodes. A supernode is a maximal set of contiguous
columns $\JS=\{j,j+1,\ldots,j+s\}$ of the $L$ factor that have the
same nonzero structure below the $(j+s)$-th row, and the lower
triangular part of $L_{\JS,\JS}$ is dense. \REV{However, this strict 
definition can produce supernodes that are either too large or too small, 
leading to memory usage, load balancing and efficiency issues.
Therefore, in our work, we relax this definition to
limit the maximal number of columns in a supernode (i.e. sets are not necessarily maximal).} 
\REV{The relaxation also allows a supernode} to include columns
for which nonzero patterns
are nearly identical to enhance the efficiency~\cite{AshcraftGrimes1989}.
\REV{This approach is also used in \superlu~\cite{LiDemmel2003}.}
Even though the nonzero pattern of the matrix can be non-symmetric, the same 
supernode partitioning is usually applied to the row partition as well, and we
assume the factorization has been computed using the structure of
$A+A^{T}$.  Then the nonzero structures of $L$ and $U$ are the transpose
of each other.  The total
number of supernodes is denoted by $\mathcal{N}$. An example of the
supernode partitioning of a structurally symmetric matrix $A$, together
with the extra fill-in in its $L$ factor ($U$ factor omitted due to
structural symmetry) are given in Fig.~\subref*{fig:matrixA.matrix}.  

Using the notation of supernodes, a pseudo-code for the selected inversion algorithm 
is given in Alg.~\ref{alg:selinvlu}.  The key step to gain computational
efficiency in the selected inversion algorithm is step 2, which
identifies the collection of all nonzero row and column indices
corresponding to the supernode $\KS$, denoted by $\CS$. 
All subsequent steps operate only on these nonzero rows and columns
within the sparsity pattern of the selected elements, thereby 
significantly reducing the computational cost.  
\begin{algorithm}
  \DontPrintSemicolon
  \caption{Selected inversion algorithm based on $LU$ factorization.}
  \label{alg:selinvlu}

  \KwIn{\begin{tabular}{l} (1) \begin{minipage}[t]{4.0in} The supernode partition of columns of $A$: $\{1,2,...,\mathcal{N}\}$ \end{minipage}\\
        (2) \begin{minipage}[t]{4.0in} A supernodal $LU$ factorization of $A$ with (unnormalized) $LU$ factors $L$ and $U$.  \end{minipage}
        \end{tabular}
        }

   \KwOut{\begin{minipage}[t]{4.0in} Selected elements of $A^{-1}$, i.e. $A^{-1}_{\IS,\JS}$ such
             that $L_{\IS,\JS}$ is not an empty block. \end{minipage}
        } 

	\For{$\KS = \REV{\mathcal{N}, \mathcal{N}-1}, ..., 1$}{
    \REV{\lnl{alg1.step0}} Find the collection of indices\;
  	$\CS=\{\IS~|~\IS>\KS,L_{\IS,\KS}\mbox{ is a
	nonzero block}\}\cup \{\JS~|~\JS>\KS, U_{\KS,\JS}\mbox{ is a
	nonzero block}\}$\;
    \REV{\lnl{alg1.step1}} $\hat{L}_{\CS,\KS}\gets L_{\CS,\KS} \REV{(L_{\KS,\KS})^{-1}},
    \hat{U}_{\KS,\CS}\gets \REV{(U_{\KS,\KS})^{-1}} U_{\KS,\CS}$\; 
  }

	\For{$\KS = \REV{\mathcal{N}, \mathcal{N}-1}, ..., 1$}{
    Find the collection of indices\;
  	$\CS=\{\IS~|~\IS>\KS,L_{\IS,\KS}\mbox{ is a
	nonzero block}\}\cup \{\JS~|~\JS>\KS, U_{\KS,\JS}\mbox{ is a
	nonzero block}\}$\;
	  \REV{\lnl{alg1.step2}} Calculate $A^{-1}_{\CS,\KS} \gets -A^{-1}_{\CS,\CS}
	  \hat{L}_{\CS,\KS}$\; 
	  \REV{\lnl{alg1.step3}} Calculate $A^{-1}_{\KS,\KS} \gets U_{\KS,\KS}^{-1}
  	L_{\KS,\KS}^{-1} - \hat{U}_{\KS,\CS} A^{-1}_{\CS,\KS}$\; 


	  \REV{\lnl{alg1.step4}} Calculate $A^{-1}_{\KS,\CS} \gets - \hat{U}_{\KS,\CS}
	  A^{-1}_{\CS,\CS}$\;
  }
\end{algorithm}

It should be noted that if $A$ is a sparse symmetric matrix, the
normalized $LU$ factors satisfy the relation 
\begin{equation}
	\hat{U}_{\CS,\KS} = \hat{L}_{\KS,\CS}^{T}. 
	\label{eqn:symul} 
\end{equation} 
\REV{To simplify the implementation, the entire diagonal block 
$A^{-1}_{\KS,\KS}$ is computed even though it is symmetric.
Due to roundoff error, the numerical update in \REV{step~\ref{alg1.step3}} may not
preserve this symmetry in finite precision.  Our numerical results
indicate that the loss of symmetry may accumulate, especially for
ill-conditioned matrices.  To reduce such error for symmetric matrices,
we can simply symmetrize the diagonal block of $A^{-1}$ by performing 
$A^{-1}_{\KS,\KS} \gets \frac12 \left( A^{-1}_{\KS,\KS} + 
A^{-T}_{\KS,\KS}\right)$ after step~\ref{alg1.step3} for each $\KS$.}

Furthermore, it should be noted that in the symmetric case, 
Eq.~\eqref{eqn:Ainv2by2normal} can be simplified using 
an $LDL^{T}$ factorization which is more efficient than
an $LU$ factorization, where $L$ is a unit lower triangular matrix, and
$D$ is a block diagonal matrix consisting of $1\times 1$ or $2\times 2$
blocks.  The simplification of the selected inversion algorithm with
an $LDL^{T}$ factorization can be found in~[Lin et al. 2011b].  For
symmetric matrices, $A^{-1}_{\KS,\CS}$ \REV{(step~\ref{alg1.step4})} is readily
obtained as the transpose of $A^{-1}_{\CS,\KS}$ without extra
computation.

\subsection{Elimination tree} \label{subsec:etree}
Both the factorization and the selected inversion can be
conveniently described in terms of traversals of an
{\em elimination tree}~\cite{Liu1990}.  Each node of the tree 
corresponds to a supernode of $A$.  A node $\RS$ is the parent 
of a node $\KS$ if and only if 
\begin{equation}
	\RS = \min\left\{ \IS > \JS ~ | ~ L_{\IS,\JS} \mbox{ is a nonzero block}
	\right\}.
	\label{}
\end{equation}
An example of the elimination tree corresponding to the matrix in
Fig.~\subref*{fig:matrixA.matrix} is given in Fig.~\subref*{fig:matrixA.etree}.

In the factorization procedure, the traversal of the elimination 
tree is a bottom-up process that starts from the leaves of the tree.
A parent supernode cannot be factored until the supernodes associated
with all its children in the tree have been factored. This type of
task dependency also determines the amount of concurrency that 
can be exploited to speed up the factorization on a parallel 
computer.

In the selected inversion procedure, the traversal of the elimination tree
is a top-down process that starts from the root of the tree.
Computing the selected elements in the $\KS$-th supernode of $A^{-1}$ 
requires the selected elements of $A^{-1}$ already computed 
at ancestor nodes of $\KS$, but not those computed 
at its sibling nodes and their descendants. 
Consequently, the selected inversion of supernodes that belong to
two different branches of the elimination tree can be in principle carried out 
independently as long as the selected elements computed at supernodes 
above these branches have been passed to processors that are assigned 
to work on these branches.

\section{Distributed memory parallel selected inversion algorithm}\label{sec:parallelization}





In this paper we present the distributed memory \pselinv method. Our
first implementation focuses on the case of symmetric matrices.
For such matrices, the selected
inversion algorithms described in Algorithm~\ref{alg:selinvlu} requires a
sparse $LU$ or $LDL^{T}$ factorization of $A$ to be available first. In
this paper we use the \superlu software package~\cite{LiDemmel2003} to
obtain the $LU$ factorization,
which has been shown to be scalable to a large number of processors 
on distributed memory parallel machines.
The relatively simple data structure of \superlu allows easy access to
sparse $L$ and $U$ factors.  However, the main ideas we develop here 
can be combined with other sparse matrix solvers such as \mumps~\cite{mumps} and
\pardiso~\cite{pardiso} too, \REV{which provides the
$LDL^{T}$ functionality and can potentially be two times faster in
the factorization phase.}
\REV{We also note that only symmetric
permutation of the matrix $A$ is allowed, even though \superlu 
allows the column permutation to be different from the row permutation.
In this work we do not perform equilibration procedure often used in the
LU factorization to modify poorly scaled matrix elements, in order to
preserve the symmetry of the matrix.}
In the current implementation of the
\pselinv method, we explicitly take advantage of the symmetry of the
matrix, and only compute the lower triangular part of the selected
elements of $A^{-1}$.  However, our implementation is not optimal in terms of memory
allocation, in the sense that both the upper and lower triangular part
of $A^{-1}$ are stored.  
As will be seen below, such a strategy
simplifies the communication pattern and efforts for bookkeeping in the
case of 2D block cyclic data distribution, and facilitates generalizing
\pselinv to asymmetric matrices.  We also note that the sub-optimal
memory allocation is not a severe limitation of the \pselinv method. The
memory footprint of \pselinv can be further optimized for applications
that are constrained by the memory usage.  Our numerical results also
indicate that the additional memory usage by \pselinv is relatively
small compared to that used in other procedures such as the \REV{parallel
numerical factorization}. 

We use the same 2D block cyclic distribution scheme 
employed in \superlu to partition and to distribute both
the $L$ factor and the selected elements of $A^{-1}$ to 
be computed. We will review the main features of this
type of distribution in Section~\ref{subsec:datadist}. 
In the 2D block cyclic distribution scheme, each
supernode $\mathcal{K}$ is assigned to and partitioned among a subset of 
processors. However, computing the selected elements of $A^{-1}$ 
contained in the supernode $\mathcal{K}$ requires retrieving 
previously computed selected elements of $A^{-1}$ that 
belong to ancestors of $\mathcal{K}$ in the elimination tree. These selected elements 
may reside on other processors. As a result, communication 
is required to transfer data among different processors to 
complete steps 3, 4 and 5 of Alg.~\ref{alg:selinvlu} 
in each iteration.  We will discuss how this is done in 
Section~\ref{subsec:block}.
%
The key to reducing communication cost and achieving scalable performance 
is to overlap communication with computation by using asynchronous 
point-to-point MPI functions, even though some of these communication 
events are collective in nature (e.g., broadcast and reduce) within 
a communication subgroup. 



In addition to utilizing a fine grain level of parallelism in 
computing $A^{-1}$ for each supernode, we introduce a coarse grain 
level of parallelism by
exploiting the concurrency available in the elimination tree. This
amounts to executing different iterates of the {\tt for} loop in
Alg.~\ref{alg:selinvlu} in parallel.  Although the elimination
tree may exhibit many independent tasks associated with supernodes that
belong to different branches of the elimination tree, the 2D block cyclic
distribution of $L$ and $A^{-1}$ may prevent these
tasks from being performed completely simultaneously on different
processors. The key to minimizing the dependency issue is to properly
assign the order of computational tasks, and to overlap computation and communication as much as
possible. We will discuss our preliminary strategy for improving the parallel
efficiency using elimination tree in Section~\ref{subsec:treeparallel}.


%
%
%

\subsection{Distributed data layout and structure} \label{subsec:datadist}

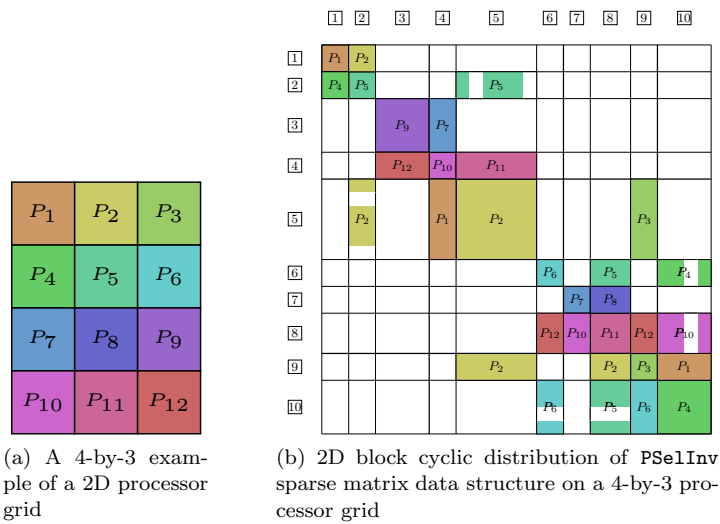
\begin{figure}
\centering

\subfloat[A 4-by-3 example of a 2D processor grid]{
\label{fig:pmatrix_grid}
\begin{adjustbox}{width=.20\linewidth}
\begin{tikzpicture}[every node/.style={font=\tiny}]

    \pgfmathtruncatemacro{\Pr}{4};
    \pgfmathtruncatemacro{\Pc}{3};

\mapblock{1}{1}{\Pr}{\Pc}
\colorproc{\curproc}{\Pr}{\Pc}
\blk{\Proc{\curproc}}{2}{2}{1}{1}{,fill=proccolor};
\mapblock{1}{2}{\Pr}{\Pc}
\colorproc{\curproc}{\Pr}{\Pc}
\blk{\Proc{\curproc}}{2}{2}{1}{3}{,fill=proccolor};
\mapblock{1}{3}{\Pr}{\Pc}
\colorproc{\curproc}{\Pr}{\Pc}
\blk{\Proc{\curproc}}{2}{2}{1}{5}{,fill=proccolor};
\mapblock{1}{4}{\Pr}{\Pc}
\colorproc{\curproc}{\Pr}{\Pc}
\blk{\Proc{\curproc}}{2}{2}{1}{7}{,fill=proccolor};

\mapblock{2}{1}{\Pr}{\Pc}
\colorproc{\curproc}{\Pr}{\Pc}
\blk{\Proc{\curproc}}{2}{2}{3}{1}{,fill=proccolor};
\mapblock{2}{2}{\Pr}{\Pc}
\colorproc{\curproc}{\Pr}{\Pc}
\blk{\Proc{\curproc}}{2}{2}{3}{3}{,fill=proccolor};
\mapblock{2}{3}{\Pr}{\Pc}
\colorproc{\curproc}{\Pr}{\Pc}
\blk{\Proc{\curproc}}{2}{2}{3}{5}{,fill=proccolor};
\mapblock{2}{4}{\Pr}{\Pc}
\colorproc{\curproc}{\Pr}{\Pc}
\blk{\Proc{\curproc}}{2}{2}{3}{7}{,fill=proccolor};

\mapblock{3}{1}{\Pr}{\Pc}
\colorproc{\curproc}{\Pr}{\Pc}
\blk{\Proc{\curproc}}{2}{2}{5}{1}{,fill=proccolor};
\mapblock{3}{2}{\Pr}{\Pc}
\colorproc{\curproc}{\Pr}{\Pc}
\blk{\Proc{\curproc}}{2}{2}{5}{3}{,fill=proccolor};
\mapblock{3}{3}{\Pr}{\Pc}
\colorproc{\curproc}{\Pr}{\Pc}
\blk{\Proc{\curproc}}{2}{2}{5}{5}{,fill=proccolor};
\mapblock{3}{4}{\Pr}{\Pc}
\colorproc{\curproc}{\Pr}{\Pc}
\blk{\Proc{\curproc}}{2}{2}{5}{7}{,fill=proccolor};

\end{tikzpicture}
\end{adjustbox}
}
~~~~~~~~
\subfloat[2D block cyclic distribution of \REV{\pselinv sparse matrix}
data structure on a 4-by-3 processor grid]{
\label{fig:pmatrix_layout}
\begin{adjustbox}{width=.45\linewidth}
\begin{tikzpicture}

    \pgfmathtruncatemacro{\Pr}{4};
    \pgfmathtruncatemacro{\Pc}{3};

\node at (-2 em, -1em) {\supernode{1}};
\node at (-2 em, -3em) {\supernode{2}};
\node at (-2 em, -6em) {\supernode{3}};
\node at (-2 em, -9em) {\supernode{4}};
\node at (-2 em, -13em) {\supernode{5}};
\node at (-2 em, -17em) {\supernode{6}};
\node at (-2 em, -19em) {\supernode{7}};
\node at (-2 em, -21.5em) {\supernode{8}};
\node at (-2 em, -24em) {\supernode{9}};
\node at (-2 em, -27em) {\supernode{10}};

\node at (1em ,2 em) {\supernode{1}};
\node at (3em ,2 em) {\supernode{2}};
\node at (6em ,2 em) {\supernode{3}};
\node at (9em ,2 em) {\supernode{4}};
\node at (13em,2 em) {\supernode{5}};
\node at (17em,2 em) {\supernode{6}};
\node at (19em,2 em) {\supernode{7}};
\node at (21.5em,2 em) {\supernode{8}};
\node at (24em,2 em) {\supernode{9}};
\node at (27em,2 em) {\supernode{10}};
\mdblk{}{1}{2}{1}{\Pr}{\Pc};
\mdblk{}{2}{2}{3}{\Pr}{\Pc};
\mdblk{}{3}{4}{6}{\Pr}{\Pc};
\mdblk{}{4}{2}{9}{\Pr}{\Pc};
\mdblk{}{5}{6}{13}{\Pr}{\Pc};
\mdblk{}{6}{2}{17}{\Pr}{\Pc};
\mdblk{}{7}{2}{19}{\Pr}{\Pc};
\mdblk{}{8}{3}{21.5}{\Pr}{\Pc};
\mdblk{}{9}{2}{24}{\Pr}{\Pc};
\mdblk{}{10}{4}{27}{\Pr}{\Pc};

\draw (0em,0em) --   (29em ,0em);
\draw (0em,-2em) --  (29em ,-2em);
\draw (0em,-4em) --  (29em ,-4em);
\draw (0em,-8em) --  (29em ,-8em);
\draw[thick] (0em,-10em) -- (29em ,-10em);
\draw (0em,-16em) -- (29em ,-16em);
\draw (0em,-18em) -- (29em ,-18em);
\draw (0em,-20em) -- (29em ,-20em);
\draw[thick] (0em,-23em) -- (29em ,-23em);
\draw (0em,-25em) -- (29em ,-25em);
\draw (0em,-29em) -- (29em ,-29em);

\draw (0em,0em) --  (0em ,-29em);
\draw (2em ,0em) -- (2em ,-29em);
\draw (4em ,0em) -- (4em ,-29em);
\draw[thick] (8em ,0em) -- (8em ,-29em);
\draw (10em,0em) -- (10em,-29em);
\draw (16em,0em) -- (16em,-29em);
\draw[thick] (18em,0em) -- (18em,-29em);
\draw (20em,0em) -- (20em,-29em);
\draw (23em,0em) -- (23em,-29em);
\draw[thick] (25em,0em) -- (25em,-29em);
\draw (29em,0em) -- (29em,-29em);

\mlublk{}{1}{2}{2}{2}{1}{3}{\Pr}{\Pc};
\mlublk{}{2}{5}{2}{6}{3}{13}{\Pr}{\Pc};

\elublkFill{}{2}{1}{3}{11.5};
\elublkFill{}{2}{1}{3}{15.5};
\draw[draw=black] ($(15.5em,-3em) + 0.5*(1em,2em)$) -- ($(15.5em,-3em) + 0.5*(1em,-2em)$);
\draw[draw=black] ($(3em,-15.5em) + 0.5*(2em,-1em)$) -- ($(3em,-15.5em) + 0.5*(-2em,-1em)$);

\mlublk{}{3}{4}{4}{2}{6}{9}{\Pr}{\Pc};
\mlublk{}{4}{5}{2}{6}{9}{13}{\Pr}{\Pc};

%
\mlublk{}{5}{9}{6}{2}{13}{24}{\Pr}{\Pc};

\mlublk{}{6}{8}{2}{3}{17}{21.5}{\Pr}{\Pc};
\mlublk{}{6}{10}{2}{4}{17}{27}{\Pr}{\Pc};

\elublkFill{}{2}{1}{17}{27.5};
\node at (17em,-27em) {\Proc{6}};
\node at (27em,-17em) {\Proc{4}};

\mlublk{}{7}{8}{2}{3}{19}{21.5}{\Pr}{\Pc};

\mlublk{}{8}{9}{3}{2}{21.5}{24}{\Pr}{\Pc};
\mlublk{+}{8}{10}{3}{4}{21.5}{27}{\Pr}{\Pc};

\elublkFill{}{3}{1}{21.5}{27.5};
\node at (21.5em,-27em) {\Proc{5}};
\node at (27em,-21.5em) {\Proc{10}};
\mlublk{}{9}{10}{2}{4}{24}{27}{\Pr}{\Pc};

%
%
%

\end{tikzpicture}
\end{adjustbox}
}

\caption{Data layout of the \REV{internal sparse matrix}
 data structure used by \pselinv.}
\end{figure}



As discussed in Section~\ref{sec:selinv}, the columns of 
$A$, $L$ and $U$ are partitioned into supernodes.
Different supernodes may have different sizes.  
The same partition is 
applied to the rows of these matrices to create a 2D block partition of these 
matrices.  The submatrix blocks are mapped to processors that are 
arranged in a virtual 2D grid of dimension $\mathrm{Pr} \times
\mathrm{Pc}$ in a cyclic fashion as follows:  The $(\IS,\JS)$-th matrix block
is held by the processor labeled by 
\REV{
\begin{equation}
	\Proc{\mathrm{mod}(\IS-1,\mathrm{Pr})\times\mathrm{Pc}+\mathrm{mod}(\JS-1,\mathrm{Pc})+1}.
	\label{}
\end{equation}
}
This is called a 2D block cyclic data-to-processor mapping.  The 
mapping itself does not take the sparsity of the matrix into account.
If the $(\IS,\JS)$-th block contains only zero elements, then that block is not
stored.  It is possible that some nonzero blocks may contain several
rows of zeros.  These rows are not stored either.
As an example, a 4-by-3 grid of processors is depicted in
Fig.~\subref*{fig:pmatrix_grid}.  
The mapping between the 2D supernode partition of the matrix 
in Fig.~\subref*{fig:matrixA.matrix} and the 2D
processor grid in Fig.~\subref*{fig:pmatrix_grid} is depicted in Fig.~\subref*{fig:pmatrix_layout}. 
Each supernodal block column of \Lmat is distributed among processors
that belong to a column of the processor grid.  Each processor may own
multiple matrix blocks.  For instance, the nonzero rows in the second 
supernode are owned by processors $\Proc{2}$ and $\Proc{5}$.
More precisely, $\Proc{2}$ owns two nonzero blocks, while $\Proc{5}$ is
responsible for one block. Note that these nonzero blocks are not necessarily contiguous
in the global matrix. 
Though the nonzero structure of \Amat is not taken into account during
the distribution, it has been shown in practice that 2D layouts leads to
higher scalability for both dense~\cite{Blackford1997}\REVTWO{,} sparse
Cholesky factorizations~\cite{RothbergGupta1994} \REV{and LU
factorization~\cite{LiDemmel2003}}.



In the current implementation, \pselinv contains an interface 
that is compatible with the \superlu software package.  
In order to allow \pselinv to be easily integrated with other
$LDL^T$ or $LU$ factorization codes, we create some intermediate
sparse matrix objects to hold the distributed $L$ and $U$ factors.  Such
intermediate sparse matrix objects will be overwritten by selected
elements of $A^{-1}$ in the selected inversion process.  Each nonzero
block \REV{$L_{\IS,\JS}$} is stored as follows. Diagonal blocks
\REV{$L_{\IS,\IS}$}
are always stored as dense matrices. Nonzero entries of
\REV{$L_{\IS,\JS}$}
($\IS>\JS$) are stored contiguously as a dense matrix in a column-major
order even though row indices associated with the stored matrix elements
are not required to be contiguous.  
%
As mentioned at the beginning of Section~\ref{sec:parallelization}, our
implementation is not optimal in terms of memory allocation for
symmetric matrices, in the sense that the nonzero entries of within
\REV{$U_{\IS,\JS}$} ($\IS<\JS$) are also stored as a dense matrix in a
contiguous array in a column major order, even though the values 
of \REV{$U_{\IS,\JS}$} are identical to those of \REV{$L^T_{\JS,\IS}$} for symmetric matrices.  The nonzero column indices
associated with the nonzeros entries in \REV{$U_{\IS,\JS}$} are not required to
be \REV{contiguous} either.  We remark that for matrices with highly asymmetric 
sparsity patterns, it is more efficient to store the upper triangular 
blocks using the skyline structure shown in~\cite{LiDemmel2003}.  However, 
we choose to use a simpler data layout because it allows level-3 BLAS (GEMM) 
to be used in the selected inversion process. 

\subsection{Computing selected elements of $A^{-1}$ within each
supernode in parallel} \label{subsec:block}

In this section, we detail how \REV{steps~\ref{alg1.step1} to~\ref{alg1.step4}} in
Alg.~\ref{alg:selinvlu} can be completed in parallel.


We perform step~\ref{alg1.step1} of Alg.~\ref{alg:selinvlu} in a separate pass,
since the data communication required in this step is relatively simple.
The processor that owns the block $L_{\KS,\KS}$ broadcasts $L_{\KS,\KS}$ to
all other processors within the same column processor group owning
nonzero blocks \REV{$L_{\IS,\KS}$} in the supernode $\KS$. Each processor in that group
performs the triangular solve \REV{$\hat{L}_{\IS,\KS}\equiv L_{\IS,\KS}
(L_{\KS,\KS})^{-1}$} for each nonzero block contained in the set
$\CS$ defined in step~\ref{alg1.step0} of the algorithm.  
Because $L_{\IS,\KS}$ is not used in the subsequent steps of selected
inversion once $\hat{L}_{\IS,\KS}$ has been computed, it is overwritten 
by $\hat{L}_{\IS,\KS}$.
Since communication is limited to a processor column group
only, step~\ref{alg1.step1} can be carried out for multiple supernodes at the same
time.


A more complicated communication pattern is required to complete step~\ref{alg1.step2} in
parallel.  Because $A^{-1}_{\CS,\CS}$ and $\hat{L}_{\CS,\KS}$ are generally
owned by different processor groups, there are two possible ways to
carry out the multiplication of $A^{-1}_{\CS,\CS}$ with $\hat{L}_{\CS,\KS}$. 
The first approach is to send blocks of $\hat{L}_{\CS,\KS}$ to
processors that own the {\em matching} blocks of $A^{-1}_{\CS,\CS}$, so
that matrix-matrix multiplication can be performed on processors owning
$A^{-1}_{\CS,\CS}$. The second approach is to send data in the opposite
direction, \ie one can send blocks of $A^{-1}_{\CS,\CS}$ to the {\em matching}
blocks of $\hat{L}_{\CS,\KS}$, so that matrix-matrix multiplication can be 
performed on processors owning $\hat{L}_{\CS,\KS}$.  

In order to compare the cost of these two approaches, 
let us first consider the case in which all blocks $\hat{L}_{\IS,\KS}$ 
with $\IS\ge \KS$ are dense matrix blocks \REV{of equal size $m_b \times
m_b$}, and $\CS=\{\IS \vert \IS\ge \KS\}$.
This is approximately the case when $\KS$ is near the root of the
elimination tree.  \REV{We assume that there are $\sqrt{P}\times \sqrt{P}$
processors, thus the size of the set $\CS$ is $m_{b}\sqrt{P}$.  
}  We also
assume that the matrix blocks in $\hat{L}_{\CS,\KS}$ are distributed
among $\sqrt{P}$ processors within the same column group, and each
processor in this processor column also holds a dense block
$\hat{L}_{\IS,\KS}$ of size $m_{b}\times m_{b}$.  In the first approach,
the computation is
performed in parallel on $P$ processors, and the computational cost on
each processor is $\mathcal{O}(m_{b}^3)$.  In the second approach, the
computation is performed in parallel on $\sqrt{P}$ processors only, and
the computational cost on these processors is
$\mathcal{O}(m_b^3\sqrt{P})$.  All other processors are idle in the
computational step, and this leads to severe load imbalance when $P$ is
large.



We implemented the first approach in \pselinv.
This requires sending the $\hat{L}_{\IS,\KS}$ block 
from a particular processor to all processors within the same column
group of processors among which $A^{-1}_{\CS,\IS}$ is distributed.
However, since the processor owning $\hat{L}_{\IS,\KS}$ is 
generally not in the same processor communication group 
that owns $A^{-1}_{\CS,\IS}$, sending $\hat{L}_{\IS,\KS}$ to processors 
that hold the distributed blocks of $A^{-1}_{\CS,\IS}$ cannot be done 
by a single broadcast.
\REV{Indeed, for this to be possible, communication groups 
(i.e. MPI communicators) would have to be created for each and every different
sparse row/column structure. This is generally not possible as the maximum number
of allowed MPI communicators is typically much smaller than needed.}
 \REV{Therefore,} one way to complete this step of data communication 
is to use a number of point-to-point MPI sends that originate 
from the processor that owns $\hat{L}_{\IS,\KS}$ and terminate on 
the group of processors that own the nonzero blocks of $A^{-1}_{\CS,\IS}$. 
In addition to incurring higher communication latency cost, this approach 
also leads to significant bookkeeping effort in order to track
the sources and destinations of all messages for each processor.

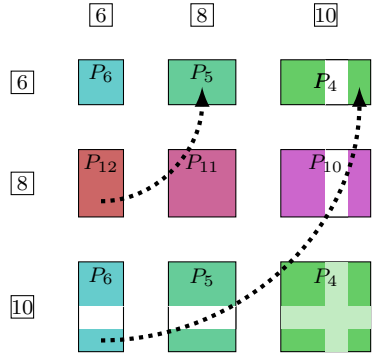
\begin{figure}[h]
\centering
\begin{adjustbox}{scale=1}
\begin{tikzpicture}

    \pgfmathtruncatemacro{\Pr}{4};
    \pgfmathtruncatemacro{\Pc}{3};

\node at (13.5em, -17em) {\supernode{6}};
\node at (13.5em, -21.5em) {\supernode{8}};
\node at (13.5em, -27em) {\supernode{10}};

\node at (17em,-14em) {\supernode{6}};
\node at (21.5em,-14em) {\supernode{8}};
\node at (27em,-14em) {\supernode{10}};


\mapblock{6}{6}{\Pr}{\Pc}
\colorproc{\curproc}{\Pr}{\Pc}
\dblknew{\begin{topbox}{2em}{2em} \vspace{2pt} \center{\Proc{\curproc}} \end{topbox}}{2}{17}{,fill=proccolor};

\mapblock{8}{8}{\Pr}{\Pc}
\colorproc{\curproc}{\Pr}{\Pc}
\dblknew{\begin{topbox}{3em}{3em} \vspace{2pt} \center{\Proc{\curproc}} \end{topbox}}{3}{21.5}{,fill=proccolor};

\mapblock{10}{10}{\Pr}{\Pc}
\colorproc{\curproc}{\Pr}{\Pc}
\dblknew{}{4}{27}{,fill=proccolor};


\mapblock{6}{8}{\Pr}{\Pc}
\colorproc{\curproc}{\Pr}{\Pc}
\blk{\begin{topbox}{2em}{3em} \vspace{2pt} \center{\Proc{\curproc}} \end{topbox}}{2}{3}{17}{21.5}{,fill=proccolor};
\mapblock{8}{6}{\Pr}{\Pc}
\colorproc{\curproc}{\Pr}{\Pc}
\blk{\begin{topbox}{3em}{2em} \vspace{2pt} \center{\Proc{\curproc}} \end{topbox}}{3}{2}{21.5}{17}{,fill=proccolor};

\mapblock{6}{10}{\Pr}{\Pc}
\colorproc{\curproc}{\Pr}{\Pc}
\blk{\begin{topbox}{2em}{4em} \vspace{2pt} \center{\Proc{\curproc}} \end{topbox}}{2}{4}{17}{27}{,fill=proccolor};
\mapblock{10}{6}{\Pr}{\Pc}
\colorproc{\curproc}{\Pr}{\Pc}
\blk{\Proc{\curproc}}{4}{2}{27}{17}{,fill=proccolor};

\mapblock{8}{10}{\Pr}{\Pc}
\colorproc{\curproc}{\Pr}{\Pc}
\blk{}{3}{4}{21.5}{27}{,fill=proccolor};
\mapblock{10}{8}{\Pr}{\Pc}
\colorproc{\curproc}{\Pr}{\Pc}
\blk{}{4}{3}{27}{21.5}{,fill=proccolor};

\mapblock{10}{10}{\Pr}{\Pc}
\colorproc{\curproc}{\Pr}{\Pc}
\blk{}{4}{1}{27}{27.5}{,fill=proccolor!40,draw=none};
\blk{}{1}{4}{27.5}{27}{,fill=proccolor!40,draw=none};

\elublkFill{}{2}{1}{17}{27.5};

\elublkFill{}{3}{1}{21.5}{27.5};
\mapblock{10}{6}{\Pr}{\Pc}
\node at (27em,-17em) {\Proc{\curproc}};

\mapblock{8}{10}{\Pr}{\Pc}
\node at (21.5em,-25.7em) {\Proc{\curproc}};
\mapblock{10}{8}{\Pr}{\Pc}
\node at (27em,-20.6em) {\Proc{\curproc}};

\mapblock{10}{10}{\Pr}{\Pc}
\node at (27em,-25.7em) {\Proc{\curproc}};

\path (17em, -22.3em) edge [dotted,ultra thick, -latex,out=0, in=-90] (21.5em, -17.2em);
\path (17em, -28.5em) edge [dotted,ultra thick, -latex,out=0, in=-90] (28.5em, -17.2em);

\end{tikzpicture}
\end{adjustbox}

\caption{Processors holding $\hat{L}_{8,6}$ and $\hat{L}_{10,6}$ send
data to processors holding the cross-diagonal blocks and overwrite
$\hat{U}_{6,8}$ and $\hat{U}_{6,10}$.  }

	\label{fig:sendtodiagonal}

\end{figure}

In our current implementation, we simplify the data communication
pattern by storing both $\hat{L}_{\IS,\KS}$ and $\hat{U}_{\KS,\IS}$ even
when $A$ is symmetric. We acknowledge that such implementation is not
optimal in terms of memory allocation, and can be improved if
applications are constrained by memory usage for symmetric matrices.
As soon as $\hat{L}_{\IS,\KS}$ becomes available as illustrated above,
we send the $\hat{L}_{\IS,\KS}$ block to the processor that owns
$\hat{U}_{\KS,\IS}$, and $\hat{U}_{\KS,\IS}$ is overwritten by
$\hat{L}_{\IS,\KS}^{T}$. Fig.~\ref{fig:sendtodiagonal} illustrates how
this step is carried out for a specific supernode $\KS=\supernode{6}$ of the matrix
described in Fig.~\subref*{fig:pmatrix_layout}.  Once $\hat{L}_{8,6}$ is
computed on $P_{12}$, the block is sent to $P_{5}$. The \REV{$P_{5}$}
processor then overwrites $\hat{U}_{6,8}$ by $\hat{L}_{8,6}^{T}$.
Similarly, the $\hat{L}_{10,6}$ block is computed on $P_{6}$ and sent to
$P_{4}$ on which $\hat{U}_{6,10}$ is overwritten by $\hat{L}_{10,6}$.

With $\hat{L}_{\IS,\KS}$ properly placed on the processors that are
mapped to the upper triangular part of the distributed $\hat{U}$ matrix,
step~\ref{alg1.step2} of Alg.~\ref{alg:selinvlu} can proceed as follows.
The $\hat{U}_{\KS,\IS}=\hat{L}^{T}_{\IS,\KS}$ block is first sent to all
processors within the same column processor group that owns
$\hat{U}_{\KS,\IS}$.  The matrix-matrix multiplication
$A^{-1}_{\JS,\IS}\hat{L}_{\IS,\KS}$ is then performed locally on each
processor owning $A^{-1}_{\JS,\IS}$ using the GEMM subroutine in BLAS3.
Then local matrix contributions $A^{-1}_{\JS,\IS}\hat{L}_{\IS,\KS}$ are 
reduced within each row communication groups owning $\hat{L}_{\JS,\KS}$ to
produce the $A^{-1}_{\JS,\KS}$ block in step~\ref{alg1.step2} of
Alg.~\ref{alg:selinvlu}. 

\begin{figure}
\centering
\begin{adjustbox}{scale=1}
\begin{tikzpicture}

    \pgfmathtruncatemacro{\Pr}{4};
    \pgfmathtruncatemacro{\Pc}{3};

\node at (10em, -17em) {\supernode{6}};
\node at (10em, -21.5em) {\supernode{8}};
\node at (10em, -27em) {\supernode{10}};

\node at (17em,-14em) {\supernode{6}};
\node at (21.5em,-14em) {\supernode{8}};
\node at (27em,-14em) {\supernode{10}};

\draw[fill=none,draw=none] (9em,-13em) rectangle +(22em,-18em);

\mapblock{6}{6}{\Pr}{\Pc}
\colorproc{\curproc}{\Pr}{\Pc}
\dblknew{\begin{topbox}{2em}{2em} \vspace{2pt} \center{\Proc{\curproc}} \end{topbox}}{2}{17}{,fill=proccolor};

\mapblock{8}{8}{\Pr}{\Pc}
\colorproc{\curproc}{\Pr}{\Pc}
\dblknew{\begin{topbox}{3em}{3em} \vspace{2pt} \center{\Proc{\curproc}} \end{topbox}}{3}{21.5}{,fill=proccolor};

\mapblock{10}{10}{\Pr}{\Pc}
\colorproc{\curproc}{\Pr}{\Pc}
\dblknew{}{4}{27}{,fill=proccolor};


\mapblock{6}{8}{\Pr}{\Pc}
\colorproc{\curproc}{\Pr}{\Pc}
\blk{\begin{topbox}{2em}{3em} \vspace{2pt} \center{\Proc{\curproc}} \end{topbox}}{2}{3}{17}{21.5}{,fill=proccolor};
\mapblock{8}{6}{\Pr}{\Pc}
\colorproc{\curproc}{\Pr}{\Pc}
\blk{\begin{topbox}{3em}{2em} \vspace{2pt} \center{\Proc{\curproc}} \end{topbox}}{3}{2}{21.5}{17}{,fill=proccolor};

\mapblock{6}{10}{\Pr}{\Pc}
\colorproc{\curproc}{\Pr}{\Pc}
\blk{\begin{topbox}{2em}{4em} \vspace{2pt} \center{\Proc{\curproc}} \end{topbox}}{2}{4}{17}{27}{,fill=proccolor};
\mapblock{10}{6}{\Pr}{\Pc}
\colorproc{\curproc}{\Pr}{\Pc}
\blk{\Proc{\curproc}}{4}{2}{27}{17}{,fill=proccolor};



\mapblock{8}{10}{\Pr}{\Pc}
\colorproc{\curproc}{\Pr}{\Pc}
\blk{}{3}{4}{21.5}{27}{,fill=proccolor};
\mapblock{10}{8}{\Pr}{\Pc}
\colorproc{\curproc}{\Pr}{\Pc}
\blk{}{4}{3}{27}{21.5}{,fill=proccolor};

\mapblock{10}{10}{\Pr}{\Pc}
\colorproc{\curproc}{\Pr}{\Pc}
\blk{}{4}{1}{27}{27.5}{,fill=proccolor!40,draw=none};
\blk{}{1}{4}{27.5}{27}{,fill=proccolor!40,draw=none};

\elublkFill{}{2}{1}{17}{27.5};

\elublkFill{}{3}{1}{21.5}{27.5};
\mapblock{10}{6}{\Pr}{\Pc}
\node at (27em,-17em) {\Proc{\curproc}};

\mapblock{8}{10}{\Pr}{\Pc}
\node at (21.5em,-25.7em) {\Proc{\curproc}};
\mapblock{10}{8}{\Pr}{\Pc}
\node at (27em,-20.6em) {\Proc{\curproc}};

\mapblock{10}{10}{\Pr}{\Pc}
\node at (27em,-25.7em) {\Proc{\curproc}};


\draw[dashed,ultra thick] (22.5em, -17em) -- (23.8em,-17em);
\draw[dashed,ultra thick, -latex] (23.8em, -17em) |- (22.5em,-21.5em);
\draw[dashed,ultra thick] (28.5em, -17em) -- (29.8em,-17em);
\draw[dashed,ultra thick, -latex] (29.8em, -17em) |- (28.5em,-21.5em);

\draw[dashed,ultra thick] (22.5em, -17em) -- (23.8em,-17em);
\draw[dashed,ultra thick, -latex] (23.8em, -17em) |- (22.5em,-27em);
\draw[dashed,ultra thick] (28.5em, -17em) -- (29.8em,-17em);
\draw[dashed,ultra thick, -latex] (29.8em, -17em) |- (28.5em,-27em);

\path (21.5em,-22.5em) edge [dotted,ultra thick, -latex,out=225, in=-45] (17em, -22.5em);
\path (21.5em,-28.5em) edge [dotted,ultra thick, -latex,out=225, in=-45] (17em, -28.5em);
\path (27em,-22.5em) edge [dotted,ultra thick, -latex,out=225, in=-45] (17em, -22.5em);
\path (27em,-28.5em) edge [dotted,ultra thick, -latex,out=225, in=-45] (17em, -28.5em);

\path (16.2em, -21.5em) edge [ultra thick, -latex,out=165, in=-180] (16.5em,-17em);
\path (16.2em, -27em) edge [ultra thick, -latex,out=165, in=-180] (16.5em,-17em);

\node at (29.8 em, -23em) {\circled{a}};
\node at (23.8 em, -23em) {\circled{a}};

\node at (21.5 em, -22em) {\circled{1}};
\node at (27 em, -22em) {\circled{1}};
\node at (21.5 em, -27.5em) {\circled{1}};
\node at (27 em, -27.5em) {\circled{1}};

\node at (19.5 em, -30em) {\circled{b}};
\node at (19.5 em, -24em) {\circled{b}};

\node at (17 em, -22em) {\circled{2}};
\node at (17 em, -27.5em) {\circled{2}};

\node at (13.8 em, -19em) {\circled{c}};
\node at (13.8 em, -19em) {\circled{c}};

\node at (17 em, -17.7em) {\circled{3}};
%
%


\end{tikzpicture}
\end{adjustbox}

\caption{Task parallelism and communication pattern for the supernode \supernode{6}.
There are 6 steps: \circled{a} broadcast $\hat{L}$, \circled{1} compute $A^{-1} \hat{L}$, \circled{b} reduce
$A^{-1} \hat{L}$, \circled{2} compute $\hat{L}^T A^{-1} \hat{L}$, \circled{c} reduce
$\hat{L}^T A^{-1} \hat{L}$ and
\circled{3} update $A^{-1}$. }

\label{fig:step3}

\end{figure}
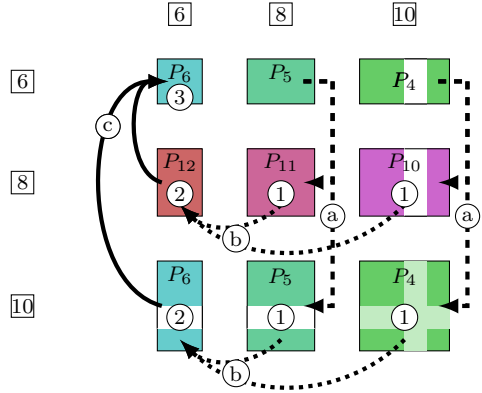

Fig.~\ref{fig:step3} illustrates how this step is completed for a
specific supernode $\KS=\supernode{6}$, for the matrix
depicted in Fig.~\subref*{fig:pmatrix_layout}.  We use circled letters
$\circled{a}, \circled{b}, \circled{c}$ to label 
communication events, and circled numbers
$\circled{1}, \circled{2}, \circled{3}$ to label computational events. 
We can see from this figure that $\hat{U}_{6,8}=\hat{L}_{8,6}^{T}$ is sent by $P_5$ 
to all processors within the same column processor group to which $P_5$ 
belongs. This group include both $P_{5}$ and $P_{11}$.  
Similarly $\hat{L}_{10,6}$ is broadcast from $P_4$ to all
other processors within the same column group to which $P_4$ belongs.
Local matrix matrix multiplications are then performed on 
$P_{11}$, $P_{10}$, $P_{4}$ and $P_5$ simultaneously. The distributed 
products are then reduced onto $P_{12}$ and $P_{5}$ within the row 
processor groups they belong to respectively.
After this step, $A^{-1}_{8,6}$ and $A^{-1}_{10,6}$ become available on
$P_{12}$ and $P_{6}$ respectively.

%
%
%
%

Upon the completion of step~\ref{alg1.step2},  the matrix product 
$\hat{U}_{\KS,\JS} A_{\JS,\KS}^{-1}\equiv\hat{L}^{T}_{\JS,\KS}
A_{\JS,\KS}^{-1}$ is first computed locally on the processor holding
$\hat{L}_{\JS,\KS}$, and then reduced to the processor that owns
the diagonal block $L_{\KS,\KS}$ within the column processor group that
the supernode $\KS$ is mapped to. The
sum of the distributed matrix product $-\hat{L}^{T}_{\JS,\KS}
A_{\JS,\KS}^{-1}$ is then added to \REV{$(U_{\KS,\KS})^{-1}
(L_{\KS,\KS})^{-1}$} computed on the processor holding $L_{\KS,\KS}$. This
completes step~\ref{alg1.step3} of Alg.~\ref{alg:selinvlu}.
As an example we use again Fig.~\ref{fig:step3} for $\KS=\supernode{6}$. $\hat{L}^{T}_{8,6}
A^{-1}_{8,6}$ is computed on $P_{12}$ and sent to $P_{6}$.  Similarly
$\hat{L}^{T}_{10,6} A^{-1}_{10,6}$ is computed on $P_{6}$.  Since both
$\hat{L}^{T}_{10,6}$ and $L_{6,6}$ are held by $P_{6}$, no further data
communication is necessary.  Finally $P_{6}$ updates $A^{-1}_{6,6}$.
As we discussed in Section~\ref{subsec:basic}, for symmetric
matrices, $A^{-1}_{\KS,\KS}$ should be explicitly symmetrized 
\REV{to reduce the effect of rounding errors}.

Since $A$ is symmetric,  
\REV{step~\ref{alg1.step4}} of Alg.~\ref{alg:selinvlu} can be simplified as follows.
We first overwrite $\hat{L}_{\JS,\KS}$ by $A^{-1}_{\JS,\KS}$ locally on the processor
holding $\hat{L}_{\JS,\KS}$ ($\JS>\KS$).  We then send
$A^{-1}_{\JS,\KS}$ to the processor holding $\hat{U}_{\KS,\JS}$, and
overwrite $\hat{U}_{\KS,\JS}$ by $(A^{-1}_{\JS,\KS})^T$.  The data
communication pattern for this step is the same as described in
Fig.~\ref{fig:sendtodiagonal}.  After \REV{step~\ref{alg1.step4}} we move to the next
supernode $(\KS-1)$.

\subsection{Exploiting concurrency in the elimination tree}\label{subsec:treeparallel}
In this section, we discuss how to add an additional coarse-grained level 
of parallelism to the selected inversion algorithm by
exploiting task concurrency exposed by the elimination tree.

As we indicated in Section~\ref{subsec:etree}, two supernodes belonging
to two separate branches of the elimination tree can be processed
independently if the selected elements of the inverse belonging to their 
ancestors have been computed, and if these supernodes and the
ancestors they depend on are mapped onto different sets of processors.  
Although it is possible to pass the
previously computed selected elements of $A^{-1}$ from the ancestors
down to their children as we move down the elimination tree, algorithms
based on this approach (e.g., a multifrontal like
algorithm~\cite{LinYangLuEtAl2011}) would require
additional work space to hold extra copies of the selected elements.

To reduce the amount of extra work space, which can grow rapidly as 
we go down the elimination tree, we choose to allow processors assigned to each
supernode to communicate back and forth with processors assigned to 
its ancestors in the way that we described in Section~\ref{subsec:block}
to complete step~\ref{alg1.step2} of Alg.~\ref{alg:selinvlu}.


However, the drawback of this approach is that, at some point, 
two supernodes belonging
to two separate branches of the elimination tree may not be processed 
simultaneously when they need to communicate with their
common ancestors at the same time. At this point of conflict, 
only one of them should be allowed to initiate and complete the data 
communication with the common ancestor at a time.  
For example, when supernodes \supernode{2} and \supernode{4} in Fig.~\subref*{fig:matrixA.matrix} are 
being processed on different sets of processors, both of them may need 
to communicate with processors assigned to supernode \supernode{5} at the 
same time.  In this case,
the updates to be performed on these processors cannot proceed 
completely independently.  \REVTWO{On the other hand, if the set of processors assigned to
update two supernodes are completely different, then at least some of
the updates  can be computed simultaneously.}


In order to exploit the type of concurrency discussed above, which occurs
at the {\tt for} loop level in Alg.~\ref{alg:selinvlu}, we create a
basic parallel task scheduler to launch different iterates of the 
{\tt for} loop in a certain order. This order is defined by
a priority list $S$, which is indexed by integer priority numbers 
ranging from 1 to $n_s$, where
$n_s$ is bounded from above by the depth of the elimination tree.
The task performed in each iteration of the {\tt for} loop is assigned 
a priority number $\sigma(\mathcal{I})$.  The lower the number, the higher the 
priority of the task, hence the sooner it is scheduled. 
The supernode $\mathcal{N}$ associated with the root of the elimination tree 
clearly has to be processed first. 
If the $k$-th element of $S$ contains multiple supernodes or tasks 
whose priority numbers are $k$,  the order in which these tasks are 
completed can be arbitrary.  
A recipe for assigning priority number of different tasks (or equivalently,
supernodes) is shown in Alg.~\ref{alg:stagemap}.  We assume that the
elimination tree is post-ordered.

Even though we use a priority list to help launch tasks, 
we do not place extra synchronization among launched tasks other than 
requiring them to preserve data dependency. 
Tasks associated with different supernodes can be executed 
concurrently if these supernodes are on different critical paths 
of the elimination tree, and if there is no overlap among processors mapped to 
these critical paths.
In fact, if tasks associated with supernode $\mathcal{J}$ and $\mathcal{I}$ are mapped to 
different sets of processors, the task associated with the supernode $\mathcal{J}$ may
actually start before that associated with another supernode $\mathcal{I}$
\REV{even if $\sigma(\mathcal{I}) < \sigma(\mathcal{J})$,}
 i.e. even if task $\mathcal{I}$ is
scheduled ahead of task $\mathcal{J}$ according to the priority list.
When two different tasks need to communicate
with a common ancestor, the priority number associated with each task
determines which task is completed first.

\begin{algorithm}[htbp]
\DontPrintSemicolon
\caption{Assign priority numbers to supernodes and create a priority list.}
\label{alg:stagemap}
				\KwIn{
						a list of supernodes $\{\mathcal{I}\}$ and the elimination tree 
						associated with these supernodes.
        }
        \KwOut{
						an array $\sigma$, $\sigma(\mathcal{I})$ gives the priority number
						of the task associated with supernode $\mathcal{I}$; an array 
						$S$ of $n_s$ supernode lists, $S(i)$ gives a set of supernodes
						with priority number $i$, $1 \leq i \leq n_s$.
        }

			  $\sigma(\mathcal{N}) = 1$\;
				$S(1) = \{\mathcal{N}\}$\;
				\For{$\mathcal{I} = \mathcal{N}-1$ down to $1$ }{
            $\sigma(\mathcal{I}) = \sigma(\textnormal{parent}(\mathcal{I})) +1$\;
            $S(\sigma(\mathcal{I})) = S(\sigma(\mathcal{I})) \cup \{ \mathcal{I} \}$\;
        }
\end{algorithm}

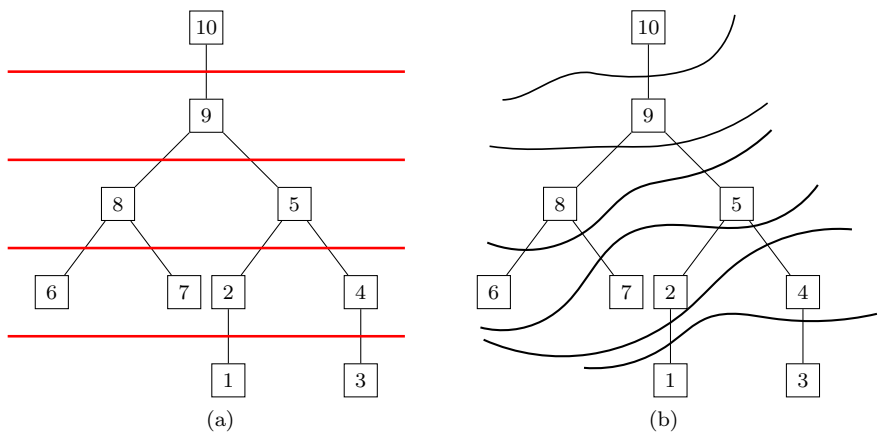
\begin{figure}[htbp]
\centering

\subfloat[]{
\label{fig:etreepart.a}
\begin{adjustbox}{width=.45\linewidth}
\begin{tikzpicture}

\node[treenode] (10) at (0,-4em) {10};
\node[treenode] (9) at (0,-8em) {9};

\node[treenode] (8) at (-4 em,-12em) {8};
\node[treenode] (5) at (4 em,-12em) {5};

\node[treenode] (6) at (-7 em,-16em) {6};
\node[treenode] (7) at (-1 em,-16em) {7};
\node[treenode] (2) at (1 em,-16em) {2};
\node[treenode] (4) at (7 em,-16em) {4};

\node[treenode] (1) at (1 em,-20em) {1};
\node[treenode] (3) at (7 em,-20em) {3};

\draw (10) -- (9);
\draw (9) -- (5);
\draw (9) -- (8);
\draw (5) -- (2);
\draw (5) -- (4);
\draw (8) -- (6);
\draw (8) -- (7);

\draw (2) -- (1);
\draw (4) -- (3);

\path[draw = red, line width = 1pt] (-9 em,-6em) --  (9 em, -6em);
\path[draw = red, line width = 1pt] (-9 em,-10em) -- (9 em, -10em);
\path[draw = red, line width = 1pt] (-9 em,-14em) -- (9 em, -14em);
\path[draw = red, line width = 1pt] (-9 em,-18em) -- (9 em, -18em);

\begin{scope}[xscale = 0.1, yscale=0.1]
\begin{scope}[cm={{0.25,0.0,0.0,0.25,(-8,0)}}]
  \path[draw=none,line join=miter,line cap=butt,miter limit=10.00,line
    width=0.640pt] (-46.9784,-87.0916) .. controls (-31.5954,-87.0613) and
    (-18.0631,-70.6571) .. (0.4720,-72.3239) .. controls (19.2815,-76.4612) and
    (53.4532,-75.5391) .. (65.9636,-64.5493) .. controls (72.5469,-58.6059) and
    (77.1635,-50.1569) .. (78.8050,-41.0444);
  \path[draw=none,line join=miter,line cap=butt,miter limit=10.00,line
    width=0.640pt] (-53.9695,-112.2290) .. controls (-30.4057,-115.6817) and
    (-14.2782,-112.5363) .. (14.4254,-112.4374) .. controls (20.5092,-112.4145)
    and (26.5952,-112.6831) .. (32.6749,-112.4374) .. controls (55.4773,-111.5159)
    and (77.9825,-103.1955) .. (96.5474,-88.8229);
  \path[cm={{0.8,0.0,0.0,-0.8,(-85.472,-40.503)}},draw=none,line join=miter,line
    cap=butt,line width=0.800pt] (37.8137,154.9439) .. controls (49.7976,159.7102)
    and (63.0976,161.1163) .. (75.8132,158.9611) .. controls (88.5289,156.8059)
    and (100.6228,151.0958) .. (110.3668,142.6467) .. controls (119.0019,135.1592)
    and (125.9271,125.5097) .. (135.8833,119.8970) .. controls (146.0021,114.1928)
    and (158.0188,113.3384) .. (169.3930,110.9816) .. controls (192.1857,106.2588)
    and (213.5673,94.9201) .. (230.2638,78.7017);
  \path[cm={{0.8,0.0,0.0,-0.8,(-85.472,-40.503)}},draw=none,line join=miter,line
    cap=butt,line width=0.800pt] (33.2023,212.4330) .. controls (47.1785,216.0514)
    and (62.4578,214.3785) .. (75.3199,207.8216) .. controls (81.4823,204.6801)
    and (87.0828,200.4675) .. (91.9211,195.5244) .. controls (99.1364,188.1529)
    and (104.5836,179.2781) .. (110.5986,170.8983) .. controls (116.6135,162.5184)
    and (123.4674,154.3886) .. (132.5016,149.4102) .. controls (141.8634,144.2513)
    and (152.8630,142.8813) .. (163.5519,142.9542) .. controls (181.5299,143.0768)
    and (199.6436,146.9814) .. (217.3518,143.8765) .. controls (234.9894,140.7839)
    and (251.2577,130.5032) .. (261.6215,115.9005);
  \path[cm={{0.8,0.0,0.0,-0.8,(-85.472,-40.503)}},draw=none,line join=miter,line
    cap=butt,line width=0.800pt] (35.3543,220.7335) .. controls (57.6548,230.6801)
    and (82.7949,234.1606) .. (106.9589,230.6468) .. controls (131.1229,227.1330)
    and (154.2306,216.6366) .. (172.7747,200.7507) .. controls (183.9002,191.2200)
    and (193.3677,179.8485) .. (204.7473,170.6227) .. controls (226.9393,152.6311)
    and (256.2110,143.6244) .. (284.6786,146.0285);
  \path[cm={{0.8,0.0,0.0,-0.8,(-85.472,-40.503)}},draw=none,line join=miter,line
    cap=butt,line width=0.800pt] (103.2959,239.7941) .. controls
    (123.7712,241.4607) and (144.7767,235.5703) .. (161.3999,223.5004) .. controls
    (169.7812,217.4148) and (177.2602,209.7236) .. (186.9164,205.9770) .. controls
    (193.7836,203.3125) and (201.3280,202.8425) .. (208.6771,203.3402) .. controls
    (216.0262,203.8380) and (223.2759,205.2669) .. (230.5712,206.2844) .. controls
    (254.1681,209.5755) and (278.3630,208.5281) .. (301.5872,203.2101);
\end{scope}
\end{scope}

\end{tikzpicture}

\end{adjustbox}
}
\subfloat[]{
\label{fig:etreepart.b}
\begin{adjustbox}{width=.45\linewidth}
\begin{tikzpicture}

\node[treenode] (10) at (0,-4em) {10};
\node[treenode] (9) at (0,-8em) {9};

\node[treenode] (8) at (-4 em,-12em) {8};
\node[treenode] (5) at (4 em,-12em) {5};

\node[treenode] (6) at (-7 em,-16em) {6};
\node[treenode] (7) at (-1 em,-16em) {7};
\node[treenode] (2) at (1 em,-16em) {2};
\node[treenode] (4) at (7 em,-16em) {4};

\node[treenode] (1) at (1 em,-20em) {1};
\node[treenode] (3) at (7 em,-20em) {3};

\draw (10) -- (9);
\draw (9) -- (5);
\draw (9) -- (8);
\draw (5) -- (2);
\draw (5) -- (4);
\draw (8) -- (6);
\draw (8) -- (7);

\draw (2) -- (1);
\draw (4) -- (3);

\path[draw = none, line width = 1pt] (-9 em,-6em) --  (9 em, -6em);
\path[draw = none, line width = 1pt] (-9 em,-10em) -- (9 em, -10em);
\path[draw = none, line width = 1pt] (-9 em,-14em) -- (9 em, -14em);
\path[draw = none, line width = 1pt] (-9 em,-18em) -- (9 em, -18em);

\begin{scope}[xscale = 0.1, yscale=0.1]
\begin{scope}[cm={{0.25,0.0,0.0,0.25,(-8,0)}}]
  \path[draw=black,line join=miter,line cap=butt,miter limit=10.00,line
    width=0.640pt] (-46.9784,-87.0916) .. controls (-31.5954,-87.0613) and
    (-18.0631,-70.6571) .. (0.4720,-72.3239) .. controls (19.2815,-76.4612) and
    (53.4532,-75.5391) .. (65.9636,-64.5493) .. controls (72.5469,-58.6059) and
    (77.1635,-50.1569) .. (78.8050,-41.0444);
  \path[draw=black,line join=miter,line cap=butt,miter limit=10.00,line
    width=0.640pt] (-53.9695,-112.2290) .. controls (-30.4057,-115.6817) and
    (-14.2782,-112.5363) .. (14.4254,-112.4374) .. controls (20.5092,-112.4145)
    and (26.5952,-112.6831) .. (32.6749,-112.4374) .. controls (55.4773,-111.5159)
    and (77.9825,-103.1955) .. (96.5474,-88.8229);
  \path[cm={{0.8,0.0,0.0,-0.8,(-85.472,-40.503)}},draw=black,line join=miter,line
    cap=butt,line width=0.800pt] (37.8137,154.9439) .. controls (49.7976,159.7102)
    and (63.0976,161.1163) .. (75.8132,158.9611) .. controls (88.5289,156.8059)
    and (100.6228,151.0958) .. (110.3668,142.6467) .. controls (119.0019,135.1592)
    and (125.9271,125.5097) .. (135.8833,119.8970) .. controls (146.0021,114.1928)
    and (158.0188,113.3384) .. (169.3930,110.9816) .. controls (192.1857,106.2588)
    and (213.5673,94.9201) .. (230.2638,78.7017);
  \path[cm={{0.8,0.0,0.0,-0.8,(-85.472,-40.503)}},draw=black,line join=miter,line
    cap=butt,line width=0.800pt] (33.2023,212.4330) .. controls (47.1785,216.0514)
    and (62.4578,214.3785) .. (75.3199,207.8216) .. controls (81.4823,204.6801)
    and (87.0828,200.4675) .. (91.9211,195.5244) .. controls (99.1364,188.1529)
    and (104.5836,179.2781) .. (110.5986,170.8983) .. controls (116.6135,162.5184)
    and (123.4674,154.3886) .. (132.5016,149.4102) .. controls (141.8634,144.2513)
    and (152.8630,142.8813) .. (163.5519,142.9542) .. controls (181.5299,143.0768)
    and (199.6436,146.9814) .. (217.3518,143.8765) .. controls (234.9894,140.7839)
    and (251.2577,130.5032) .. (261.6215,115.9005);
  \path[cm={{0.8,0.0,0.0,-0.8,(-85.472,-40.503)}},draw=black,line join=miter,line
    cap=butt,line width=0.800pt] (35.3543,220.7335) .. controls (57.6548,230.6801)
    and (82.7949,234.1606) .. (106.9589,230.6468) .. controls (131.1229,227.1330)
    and (154.2306,216.6366) .. (172.7747,200.7507) .. controls (183.9002,191.2200)
    and (193.3677,179.8485) .. (204.7473,170.6227) .. controls (226.9393,152.6311)
    and (256.2110,143.6244) .. (284.6786,146.0285);
  \path[cm={{0.8,0.0,0.0,-0.8,(-85.472,-40.503)}},draw=black,line join=miter,line
    cap=butt,line width=0.800pt] (103.2959,239.7941) .. controls
    (123.7712,241.4607) and (144.7767,235.5703) .. (161.3999,223.5004) .. controls
    (169.7812,217.4148) and (177.2602,209.7236) .. (186.9164,205.9770) .. controls
    (193.7836,203.3125) and (201.3280,202.8425) .. (208.6771,203.3402) .. controls
    (216.0262,203.8380) and (223.2759,205.2669) .. (230.5712,206.2844) .. controls
    (254.1681,209.5755) and (278.3630,208.5281) .. (301.5872,203.2101);
\end{scope}
\end{scope}

\end{tikzpicture}

\end{adjustbox}
}
\caption{Elimination tree of matrix $A$ and two possible priority lists $S$.}
\label{fig:etreepart}
\end{figure}

We remark that there is some flexibility in assigning a priority number to 
each supernode and constructing the priority list $S$.
For instance, we can use the strategy given by
Alg.~\ref{alg:stagemap}, which simply defines $\sigma(\IS)$ by the
distance (in terms of the number of edges) between the supernode $\IS$ and the root of the elimination tree.
For the same elimination tree shown in Fig.~\subref*{fig:etreepart.a}, another 
possible construction of the $\sigma$ list is illustrated in
Fig.~\subref*{fig:etreepart.b}, which assigns the same $\sigma$ value to 
supernodes at different levels of the elimination tree. The latter 
construction takes into account how supernodes are distributed among
different processors as we will discuss below.

The priority list $S$ determines the order in which 
computational tasks associated with different supernodes are completed.
Because the amount of work and communication performed by each
supernodes can vary significantly, different priority lists
can lead to different overall performance.
The actual performance of parallel selected inversion depends on 
the sparsity pattern of the matrix as well as the processor grid, 
and is therefore difficult to predict \textit{a priori}. \REVTWO{We refer
readers to our report~\cite{JacquelinLinYang2014Report} for a detailed
example on the difficulty of designing an optimal task schedule.}

With the help of the priority list, we can implement the {\tt for} loop level of
parallelism in Alg.~\ref{alg:selinvlu} in a way that is described in
Alg.~\ref{alg:pselinvldlt}. \REV{To illustrate the relation clearly,
Alg.~\ref{alg:pselinvldlt} uses the same numeric ordering of the steps
as that in Alg.~\ref{alg:selinvlu}.}
Alg.~\ref{alg:pselinvldlt} also makes use of another array of lists 
{\tt procmap}.  The $\KS$-th element of {\tt procmap} 
contains the list of all processors participating in steps \REV{$\ref{alg1.step2}-\ref{alg1.step4}$} in
Alg.~\ref{alg:selinvlu}.  The communication
steps are described within parentheses.  We also remark that we do not
place MPI barriers between supernodes explicitly to exploit parallelism
among the computation for different supernodes.  For symmetric matrices,
the diagonal blocks should be symmetrized, as indicated \REVTWO{at the end of
Section~\ref{subsec:basic} for symmetric cases.}

\begin{algorithm}[htbp]
\DontPrintSemicolon

  \KwIn{ 
            \begin{tabular}{l}
            (1) 
                 \begin{minipage}[t]{4.0in}
                 The supernode partition of columns of a sparse
								 symmetric matrix $A$: $\{1,2,\ldots,\mathcal{N}\}$; a priority
                 list $\{S(k)\}$: $k=1,2,\ldots,n_{s}$;
                  \end{minipage} \\
            (2) \begin{minipage}[t]{4.0in}
                  $L$ and $U$ factors through a supernodal $LU$
                  factorization (or equivalent $LDL^{T}$ factorization)
                  of $A$; 
                  \end{minipage} \\
            (3) \begin{minipage}[t]{4.0in}
                2D processor mapping with with $P=\mathrm{Pr}\times
                \mathrm{Pc}$ processors.
                  \end{minipage}
            \end{tabular}
          }
  \KwOut{  Selected elements of $A^{-1}$.}

  [Compute the normalized factors $\hat{L}$ and
  $\hat{U}$].\;

  \For{\REVTWO{$k$ = 1, 2,..., $n_s$}}{
  \For {each supernode $\KS \in S(k)$}{
  \If{myid $\in$ {\tt procmap(\KS)}}{
  \REV{\nlset{\ref{alg1.step0}}} Find the collection of indices
  $\CS=\{\IS~|~\IS>\KS,L_{\IS,\KS}\mbox{ is a
  nonzero block}\}\cup \{\JS~|~\JS>\KS, U_{\KS,\JS}\mbox{ is a
  nonzero block}\}$\;
  \REV{\nlset{\ref{alg1.step1}}} (Broadcast \REV{$(L_{\KS,\KS})^{-1}$} to processors owning
  $L_{\IS,\KS},\IS\in\CS$)\;
  \REV{\nlset{\ref{alg1.step1}}} $\hat{L}_{\CS,\KS}\gets L_{\CS,\KS} \REV{(L_{\KS,\KS})^{-1}}$\;
  \REV{\nlset{\ref{alg1.step1}}} (Send $\hat{L}_{\IS,\KS},\IS\in\CS$ to the processor holding
  $\hat{U}_{\KS,\IS}\REVTWO{,\IS\in\CS}$ and overwrite $\hat{U}_{\KS,\IS}\REVTWO{,\IS\in\CS}$ by
  $\hat{L}^T_{\IS,\KS}\REVTWO{,\IS\in\CS}$)\;
  }
  }
  }
  [Selected inversion process].\;
  \For{$k$ = 1, 2,..., $n_s$}{
  \For{each supernode $\KS \in S(k)$}{
  \If{myid $\in$ {\tt procmap(\KS)}}{
  Find the collection of indices
  $\CS=\{\IS~|~\IS>\KS,L_{\IS,\KS}\mbox{ is a
  nonzero block}\}\cup \{\JS~|~\JS>\KS, U_{\KS,\JS}\mbox{ is a
  nonzero block}\}$\;
  \REV{\nlset{\ref{alg1.step2}}} (Broadcast $\hat{U}_{\KS,\IS}\REVTWO{,\IS\in\CS}$ to processors holding
  $A^{-1}_{\JS,\IS}\REVTWO{,\IS},\JS\in\CS$)\;
  \REV{\nlset{\ref{alg1.step2}}} For processors holding $A^{-1}_{\JS,\IS},\IS,\JS\in\CS$, compute
  locally $-A^{-1}_{\JS,\IS} \hat{L}_{\IS,\KS}\REVTWO{,\IS,\JS\in\CS}$\; 
  \REV{\nlset{\ref{alg1.step2}}} (Reduce $-A^{-1}_{\JS,\IS} \hat{L}_{\IS,\KS}\REVTWO{,\IS,\JS\in\CS}$ to processors
  holding $\hat{L}_{\JS,\KS}\REVTWO{,\JS\in\CS}$, and save the result in $A^{-1}_{\JS,\KS}\REVTWO{,\JS\in\CS}$)\;
  \REV{\nlset{\ref{alg1.step3}}} For processors holding \REV{$\hat{L}_{\IS,\KS},\IS\in\CS$}, compute
  locally $-\hat{L}_{\IS,\KS}^{T} A^{-1}_{\IS,\KS}\REVTWO{,\IS\in\CS}$\; 
  \REV{\nlset{\ref{alg1.step3}}} (Reduce $-\hat{L}_{\IS,\KS}^{T} A^{-1}_{\IS,\KS}\REVTWO{,\IS\in\CS}$ to the
  processor holding $L_{\KS,\KS}$)\;
  \REV{\nlset{\ref{alg1.step3}}} For the processor holding \REV{$(L_{\KS,\KS})^{-1}$}, update
  \REV{$A^{-1}_{\KS,\KS} \gets
  (U_{\KS,\KS})^{-1} (L_{\KS,\KS})^{-1} - \hat{L}_{\CS,\KS}^{T}
  A^{-1}_{\CS,\KS}$}\;
  For the processors holding $L^{-1}_{\KS,\KS}$, update
  $A^{-1}_{\KS,\KS} \gets \frac12 \left( A^{-1}_{\KS,\KS} + 
  A^{-T}_{\KS,\KS}\right)$\;
  \REV{\nlset{\ref{alg1.step2}}} Overwrite $\hat{L}_{\IS,\KS}\REVTWO{,\IS\in\CS}$ by $A^{-1}_{\IS,\KS}\REVTWO{,\IS\in\CS}$\;
  \REV{\nlset{\ref{alg1.step4}}} (Send $A^{-1}_{\IS,\KS},\IS\in\CS$ to the processor holding 
  $\hat{U}_{\KS,\IS}\REVTWO{,\IS\in\CS}$, and overwrite $\hat{U}_{\KS,\IS}\REVTWO{,\IS\in\CS}$ by $A^{-1}_{\IS,\KS}\REVTWO{,\IS\in\CS}$)\;
  }
}
}
\caption{The parallel selected inversion algorithm (for symmetric
matrices).}
\label{alg:pselinvldlt}
\end{algorithm}

\section{Numerical results}\label{sec:numerical}

To assess the performance of \pselinv, we conducted a 
number of computational experiments which we report in this section.

Our test problems are taken from various sources including Harwell-Boeing Test Collection
\cite{HarwellBoeing}, the University of Florida Matrix
Collection\cite{FloridaMatrix}, and matrices generated from electronic
structure software including SIESTA~\cite{SolerArtachoGaleEtAl2002} and
DGDFT~\cite{LinLuYingE2012}. 
The first two matrix collections are widely used benchmark problems for
testing sparse direct methods, while the other test problems come from
practical large scale electronic structure calculations.  The names
of these matrices as well as some of their characteristics are listed in
Tables~\ref{tab:testproblem} and \ref{tab:characteristics}.  The first
three problems in these tables come with two matrices each. One of the 
matrices, denoted by $H$, is a discretized Hamiltonian, and the other matrix
is an overlap matrix denoted by $S$. For all other problems, the overlap
matrices can be considered as the identity matrix.
All matrices are real and symmetric. In all our experiments, we compute
the selected elements of the matrix
\begin{equation}
A(z) = H - z S.
\label{eq:az}
\end{equation}
For simplicity, we choose $z=0$ for all the efficiency tests in
Section~\ref{subsec:scalability}.  \REV{For some applications,
$z$ is chosen to be a complex number with a small imaginary part
to ensure that $A(z)$ is nonsingular. This technique is often used in 
electronic structure calculation to be discussed in Section~\ref{subsec:ksdft}.}
The $LU$ factorization is performed
by using the \superlu software package.  
\superlu does not use dynamic pivoting strategies and our matrices are permuted
without taking into account the values of matrix entries. Consequently
the efficiency of both \superlu and \pselinv is independent of the
choice of $z$.  The lack of dynamic pivoting strategies may impact the
accuracy of \pselinv for highly indefinite and nearly singular systems.
We study the accuracy for different choices of complex shifts $z$ in
Section~\ref{subsec:accuracy}.  All the timing results reported
are performed in complex arithmetic computation.

\begin{table}
  \centering
  \begin{tabular}{l|l}
    \toprule
    Problem & Description \\
    \midrule
		\REV{\lucbnc} & \REV{Electronic structure theory, C-BN sheet with
		2532 atoms} \\ 
		SIESTA\_C\_BN\_2x2 & Electronic structure theory, C-BN sheet with
		10128 atoms \\ 
		SIESTA\_C\_BN\_4x2 & Electronic structure theory, C-BN sheet with
		20256 atoms \\ 
    \REV{\ludna} & \REV{Electronic structure theory, DNA molecule with
    11440 atoms} \\
    DNA\_715\_64cell & Electronic structure theory, DNA molecule with
    45760 atoms \\
		DG\_Graphene\_2048 & Electronic structure theory, graphene with 2048
		atoms \\
		DG\_Graphene\_8192 & Electronic structure theory, graphene with 8192
		atoms \\
    pwtk & Pressurized wind tunnel, stiffness matrix. \\
    parabolic\_fem & Diffusion-convection reaction, constant homogeneous
    diffusion. \\
    ecology2 & Circuitscape: circuit theory applied to animal/gene flow,
    B. McRae, UCSB. \\
    audikw\_1 & Automotive crankshaft model with over 900,000 TETRA elements, Audi, GmbH. \\
    \bottomrule
  \end{tabular}
  \caption{Description of test problems for \pselinv.}
  \label{tab:testproblem}
\end{table}

\begin{table}
  \centering
  \begin{tabular}{c|c|c|c}
    \toprule
    Problem & $n$ & $|A|$ & $|L|$ \\
    \midrule
		\REV{\lucbnc} & \REV{32,916} & \REV{23,857,418} &  \REV{269,760,112}\\
		SIESTA\_C\_BN\_2x2 & 131,664 & 95,429,672 &  1,655,233,542\\
		SIESTA\_C\_BN\_4x2 & 263,328  & 190,859,344  & 3,591,750,262 \\
    \REV{\ludna} & \REV{7,752}& \REV{2,430,642}& \REV{9,272,160}\\
    DNA\_715\_64cell & 459,712& 224,055,744& 866,511,698\\
		DG\_Graphene\_2048 & 82,944  & 87,340,032  &545,245,344 \\
		DG\_Graphene\_8192 & 331,776  & 349,360,128  & 2,973,952,468\\
    pwtk & 217,918 & 5,926,171 & 104,644,472  \\
    parabolic\_fem & 525,825 & 3,674,625 & 58,028,731\\
    ecology2 & 999,999 & 2,997,995 & 91,073,583 \\
    audikw\_1 & 943,695& 77,651,847& 2,500,489,909\\
    \bottomrule
  \end{tabular}
  \caption{The dimension $n$, the number of nonzeros $|A|$, and the
  number of nonzeros of the Cholesky factor $|L|$ of the test problems.}
  \label{tab:characteristics}
\end{table}

\begin{figure}[htbp]
\centering
\begin{adjustbox}{width=.8\linewidth}
\input{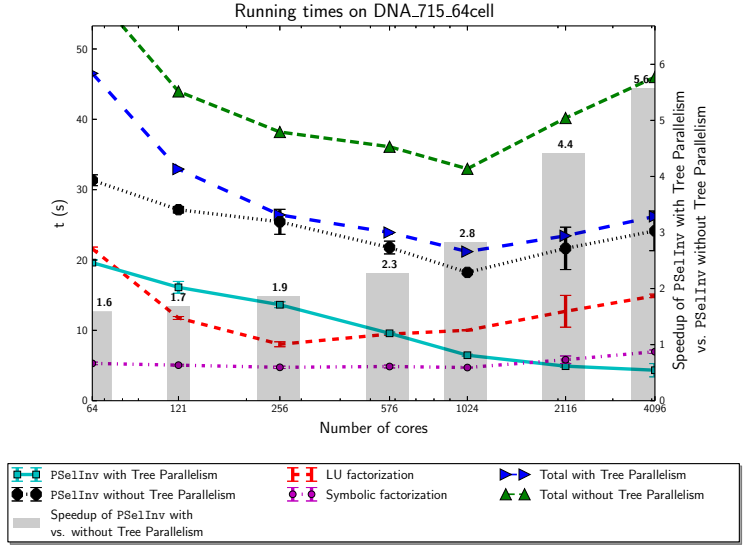}
\end{adjustbox}
\caption{The wall clock time used by three components (symbolic 
factorization, numerical LU factorization and selected inversion) 
of two versions of \pselinv with respect \REV{to} the number of \REVTWO{cores} 
used in the computation for the DNA\_715\_64cell matrix. In one version, 
we do not take advantage of the concurrency exposed by the elimination 
tree, whereas in the other version we do make use of this tree level of 
parallelism. The height of each bar in the figure indicates the ratio of 
wall clocked time measured for the former over that measured for the 
latter. }
\label{fig:DNA64_treeimpact}
\end{figure}

In all of our experiments, we used the NERSC Edison platform with Cray
XC30 nodes. Each node has 24 \REVTWO{cores} partitioned among two
Intel Ivy Bridge processors.  Each 12-core processor runs at 2.4GHz. A
single node has 64GB of memory, providing more than 2.6 GB of memory per
core. \REVTWO{We used one MPI process per core, and refer to core to denote
an MPI process.} 


\subsection{Parallelization scalability}\label{subsec:scalability}
\begin{figure}[htbp]
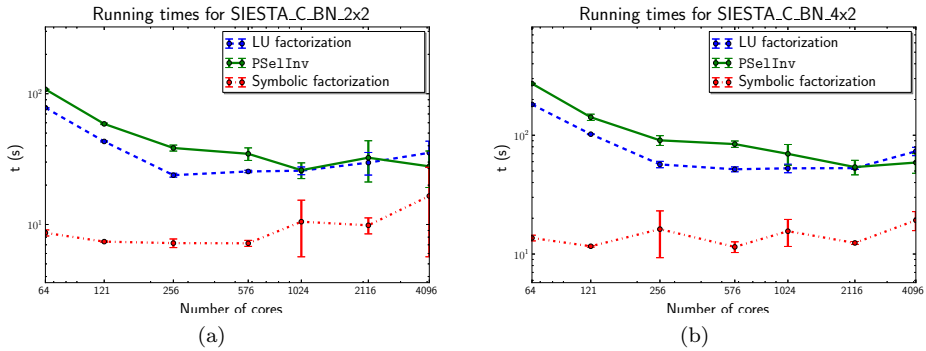

\centering

\subfloat[]{%
\centering%
\begin{adjustbox}{width=.48\linewidth}%
\input{ss_LU_C_BN_C_2by2.pgf}%
\end{adjustbox}%
}%
\quad
\subfloat[]{
\begin{adjustbox}{width=.48\linewidth}%
\input{ss_LU_C_BN_C_4by2.pgf}%
\end{adjustbox}%
}
\caption{The strong scalability of \pselinv compared to that of LU
factorization and symbolic factorization for SIESTA matrices generated 
for \REV{two} C\_BN systems of different sizes.}\label{fig:scalability_siesta}
\end{figure}

\begin{figure}[htbp]
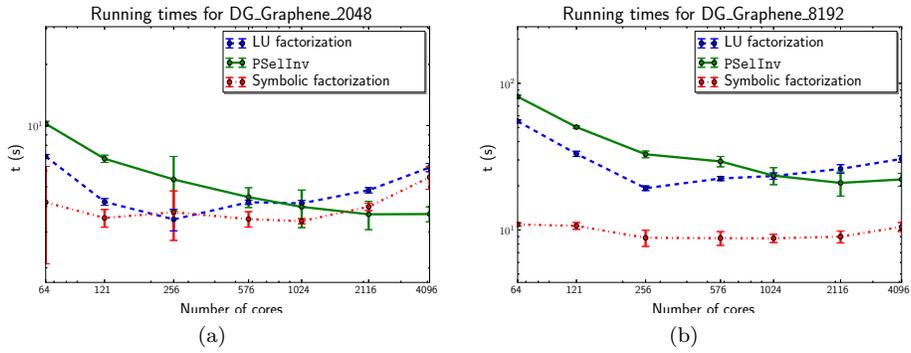

\centering
\subfloat[]{
\begin{adjustbox}{width=.48\linewidth}
\input{ss_DGGraphene_2048.pgf}
\end{adjustbox}
}
\subfloat[]{
\begin{adjustbox}{width=.48\linewidth}
\input{ss_DGGraphene_8192.pgf}
\end{adjustbox}
}

\caption{The strong scalability of \pselinv compared to that of LU
factorization and symbolic factorization for DG matrices generated 
for graphene systems of different sizes.}
\label{fig:scalability_dgdft}
\end{figure}

\begin{figure}[htbp]
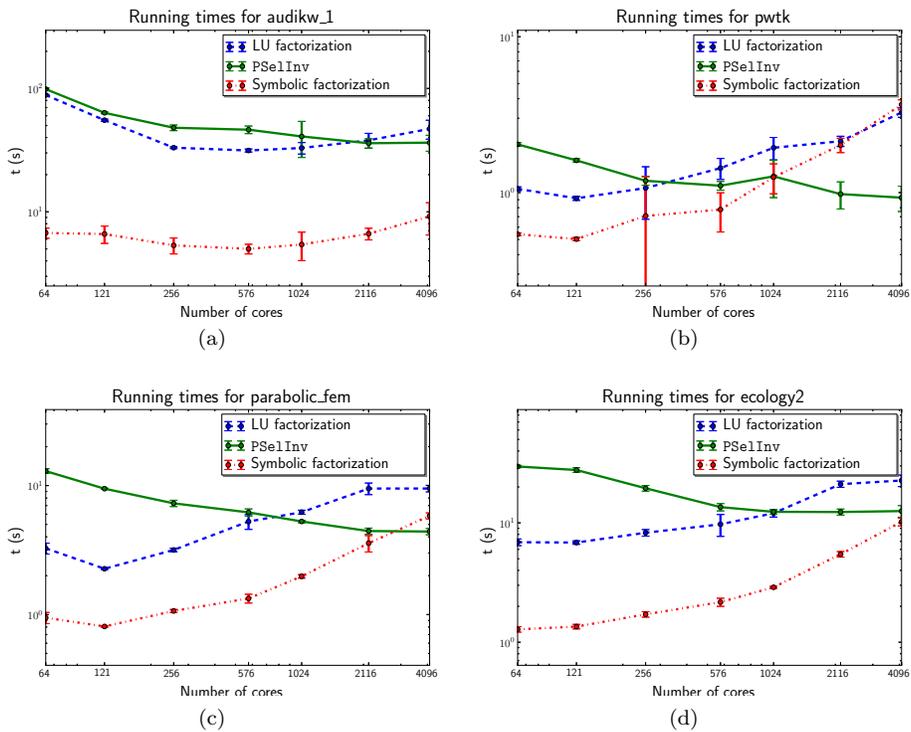

\centering

\subfloat[]{
\begin{adjustbox}{width=.48\linewidth}
\input{ss_Audikw_1.pgf}
\end{adjustbox}
}
\subfloat[]{
\begin{adjustbox}{width=.48\linewidth}
\input{ss_Pwtk.pgf}
\end{adjustbox}
}

\subfloat[]{
\begin{adjustbox}{width=.48\linewidth}
\input{ss_Parabolic_fem.pgf}
\end{adjustbox}
}
\subfloat[]{
\begin{adjustbox}{width=.48\linewidth}
\input{ss_Ecology2.pgf}
\end{adjustbox}
}
\caption{The strong scalability of \pselinv compared to that of LU
factorization and symbolic factorization for matrices from
Harwell-Boeing Test Collection and the University of Florida Matrix
Collection.}
\label{fig:scalability_ufl}
\end{figure}

\REV{Four} types of experiments were performed to measure the scalability
of \pselinv. \REV{For the first three sets of experiments,} all the timing
data points we present are averaged measurements over 10 runs, and the error bars shown in 
Figgures~\ref{fig:DNA64_treeimpact}, ~\ref{fig:scalability_siesta},
~\ref{fig:scalability_dgdft} and ~\ref{fig:scalability_ufl} 
indicate the standard deviation of the measured wall clock time.   

In order to clearly show the cost and scalability of selected inversion
itself in comparison with the symbolic and numerical LU factorizations, 
which are required for selected inversion, we time the three computational
components separately.
The symbolic factorization is performed in parallel using \ptscotch. 
It is labeled by ``symbolic factorization'' in the timing figures.  
The LU factorization is performed by using \superlu. It is labeled by 
``LU factorization''. The selected inversion itself is performed by using 
\pselinv, and labeled by ``PSelInv''.
The total time required to obtain the selected elements of the inverse 
matrix thus corresponds to the sum of all three components.

The first experiment focuses on the impact of the additional parallelism
stemming from the elimination tree as discussed in
Section~\ref{subsec:treeparallel}. \pselinv is thus tested both with and
without this additional level of parallelism. As observed in
Fig.~\ref{fig:DNA64_treeimpact}, adding the tree level parallelism 
\REV{allows the performance of \pselinv to scale to 4096 \REVTWO{cores}. 
When 4096 \REVTWO{cores} are used, adding the tree level parallelism 
leads to a 5.6 fold speedup in for the DNA\_715\_64cell 
problem.  Comparatively, the performance of the LU factorization and the 
selected inversion without tree parallelism can only scale up to 
$1024$ \REVTWO{cores}.}

The second set of experiments
(Figures.~\ref{fig:scalability_siesta}, \ref{fig:scalability_dgdft} and \ref{fig:scalability_ufl}) aims at evaluating the strong scaling of \pselinv.
In every experiment, the tree level of parallelism is enabled as it clearly 
delivers better performance.
\pselinv exhibits excellent strong scalability up to 
4,096 \REVTWO{cores}.  For
SIESTA matrices (Fig.~\ref{fig:scalability_siesta}), \pselinv is slightly 
slower than LU factorization when the number of \REVTWO{cores} is less than $1024$, and is faster than the LU factorization when more than $2116$ \REVTWO{cores} 
are used.  For DGDFT matrices (Fig.~\ref{fig:scalability_dgdft}),
\pselinv can be twice as fast as the LU factorization, and the running
time of \pselinv can be comparable to that of symbolic factorization
for the DNA\_715\_64cell matrix.  For
generic sparse matrices (Fig.~\ref{fig:scalability_ufl}) obtained from
the University of Florida Collection, \pselinv delivers excellent
performance on relatively dense matrices, such as audikw\_1 and pwtk.
We also observe that for highly sparse problems, such as ecology\_2 and parabolic\_fem,
\pselinv is relatively more costly, but the scalability of \pselinv can
still be better than that of \superlu when a large number of \REVTWO{cores} are used.


\begin{figure}[htbp]
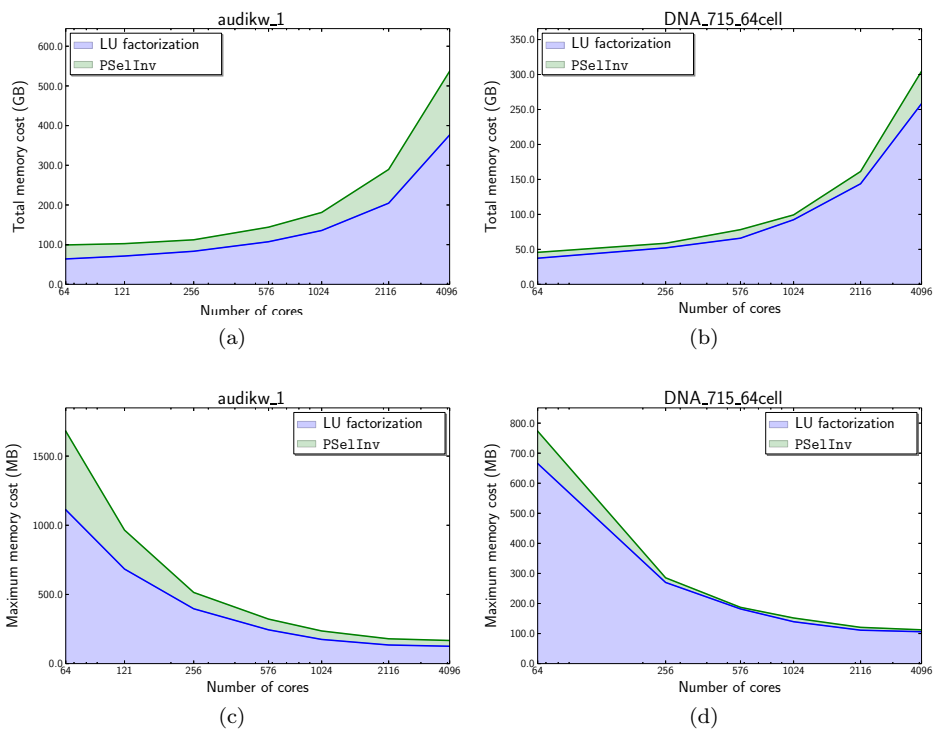

\centering

\REV{
\subfloat[]{
\label{fig:mem_footprint.Audikw_1total}
\begin{adjustbox}{width=.48\linewidth}
\input{mc_Audikw_1.pgf}
\end{adjustbox}
}
\subfloat[]{
\label{fig:mem_footprint.DNA64total}
\begin{adjustbox}{width=.48\linewidth}
\input{mc_DNA_715_64cell.pgf}
\end{adjustbox}
}

\subfloat[]{
\label{fig:mem_footprint.Audikw_1percore}
\begin{adjustbox}{width=.48\linewidth}
\input{mc2_Audikw_1.pgf}
\end{adjustbox}
}
\subfloat[]{
\label{fig:mem_footprint.DNA64percore}
\begin{adjustbox}{width=.48\linewidth}
\input{mc2_DNA_715_64cell.pgf}
\end{adjustbox}
}
}

\caption{The cumulative total memory cost and maximum memory cost among all
\REVTWO{cores} for the audikw\_1 matrix and the DNA\_715\_64cell matrix,
with shaded region indicating the additional memory cost introduced in
each step of the selected inversion process. The memory cost is measured
by memory high watermark reached at each step.}
\label{fig:mem_footprint}

\end{figure}

The third experiment focuses on the memory cost of \pselinv.
Fig.~\ref{fig:mem_footprint} shows the cumulative total
memory cost and cumulative maximum single core memory usage in the parallel selected inversion procedure for two matrices: 
audikw\_1 and DNA\_715\_64cell. 
\REV{This is measured by the memory high watermark reached after the LU
factorization step and the selected inversion step, respectively. The
memory cost for reading the input matrix is much smaller compared to the
memory cost of LU factorization or selected inversion.
The shaded areas in
Fig.~\ref{fig:mem_footprint} correspond to the additional memory
required to perform each step of the procedure.  
The memory high watermark implies that the memory cost for 
LU factorization not only includes the memory cost for the LU factors,
but also other temporary memory allocation created during the LU
factorization process.  The situation is similar for the selected
inversion process.
The memory required to store both the LU factors and the corresponding
elements in the inverse is lower than the overall memory high watermark,
which accounts for additional memory required to hold communication and
temporary buffers.  
}
We can see from Fig.~\ref{fig:mem_footprint}, most of the memory allocation
is done during the LU factorization step. 
The total memory cost of LU factorization and selected
inversion increases as the number of \REVTWO{cores} increases due to the use of
additional buffer arrays for communication and computation.

The total additional memory cost of \pselinv is $20\%\sim 60\%$ of the total 
memory required by LU factorization, which is relatively small.
The maximum memory usage per core decreases steadily as the number 
of \REVTWO{cores} increases.  For our test problems, the maximum memory cost per core of \pselinv
is around 1GB when a relatively small number (64) of \REVTWO{cores} are
used. It decreases to around 100MB when a large number (4096) of
\REVTWO{cores} are used for the same problem.

\REV{The last set of experiments compare the performance of \pselinv with 
the parallel matrix inversion method recently implemented in the \mumps
package~\cite{AmestoyDuffLExcellentEtAl2012a}.
In the MUMPS algorithm, the actual set of computed entries of $A^{-1}$
is also a \textit{superset} of the requested entries of $A^{-1}$.  This
method can be more efficient than \pselinv when a small number of entries of
$A^{-1}$ are requested.  However, when a relatively large number of
entries are to be computed such as in the computation of the selected
elements defined by~\eqref{eqn:selelem}, \pselinv can more efficiently
reuse the information shared among different entries of
$A^{-1}$.
Fig.~\ref{fig:perf_mumps} shows 
the numerical factorization and selected inversion timings
for
\mumps and \pselinv, respectively.  As all matrices are symmetric, 
MUMPS performs $LDL^T$ factorizations, which use
fewer floating point operations and less memory.
\REVTWO{We use \mumps to compute the
selected elements of $A^{-1}$ as defined in Eq.~\eqref{eqn:selelem} for
\ludna and \lucbnc, and use \mumps to only compute the diagonal entries
of $A^{-1}$ for \dgacpn. All three matrices come from electronic
structure applications, and the difference in terms of computed entries
is
determined by the different requirements in practical calculations.
\mumps contains a block size parameter that controls the number of
right-hand sides processed simultaneously. We experimented with this
parameter, and report the best results we could produce for each case
(i.e. 16 for \ludna and 256 for \lucbnc and \dgacpn)}.
Table~\ref{tab:mumpscomp} shows how the wallclock time used
by 
\pselinv compares with that used by MUMPS for \REVTWO{three} test problems
\ludna, \REVTWO{\dgacpn} and \lucbnc.

%

\begin{table}
  \centering
\REVTWO{
  \begin{tabular}{c|c|c|c}
    \toprule
    Problem & \# of \REVTWO{cores} & \pselinv & \mumps \\
    \midrule
\ludna  & 1 & 3.1  & 29.9 \\
\ludna  & 16 & 2.9 & 44.6\\
\dgacpn & 1  & 4.9 & 82.6\\
\dgacpn & 16 & 0.6 & 67.1\\
\lucbnc & 16  & 52.6 & 648.2 \\
\lucbnc & 256 & 5.6 & 293.7 \\
    \bottomrule
  \end{tabular}
  \caption{\REV{A comparison of wallclock time required by \pselinv and
\mumps for computing the selected elements defined by~\eqref{eqn:selelem}
of two matrices \ludna and \lucbnc\REVTWO{, and the diagonal elements of \dgacpn}.} }
  \label{tab:mumpscomp}
}
\end{table}

For \REVTWO{all} test problems, \pselinv is at least an order of magnitude
faster than \mumps. It is two orders of magnitude faster on \lucbnc when
using 256 \REVTWO{cores}.  \pselinv can also scale to a relatively
larger number of \REVTWO{cores}. Our numerical results indicate that for
the computation of selected elements as considered in this paper,
\pselinv is more efficient than the more general approach taken in
\mumps. 

\REVTWO{One way to understand the speedup of \pselinv over \mumps is
through the following idealized situation. Consider the computation of
the diagonal entries of a tridiagonal matrix $A$ of size $N$, which can
be computed by solving $N$ set of triangular equations of the
form~\eqref{eqn:triangular}. Even with the help of the elimination tree
and the fact that only one entry of $A^{-1}$ is needed for each
equation, the cost for solving all $N$ equations independently would be
$\Or(N^2)$. The reason for the high computational cost is that a
sigificant amount information calculated among the $N$ equations are
redundant.  The \mumps approach is to detect such redundant information
on the fly through graph based algorithms.  However, designing an
optimal algorithm which maximally reduces the amount of redundant
information for a given data-to-processor mapping is difficult.  The
resulting implementation may or may not reach the optimal complexity.
By contrast \pselinv removes all the redundant calculation \textit{by
design}, and the computational cost for a tri-diagonal matrix can be
provably reduced to $\Or(N)$~[Lin et al. 2009b].}

We should point out that \mumps can compute an arbitrary set of 
elements of $A^{-1}$, whereas the set of selected elements that
can be computed by \pselinv is more restrictive as defined in
Eq.~\eqref{eqn:selelem}.

\begin{figure}[htbp]
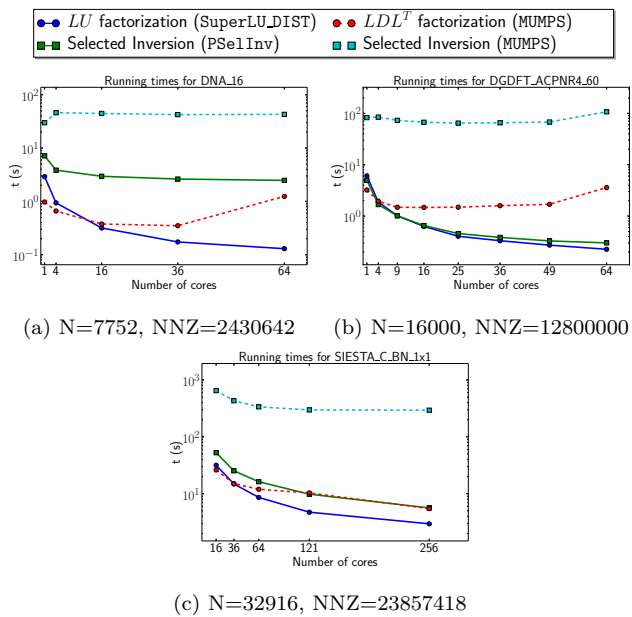

\centering
\REVTWO{
\begin{adjustbox}{width=.6\linewidth}
\begingroup%
\makeatletter%
\begin{tikzpicture}%
\pgfpathrectangle{\pgfqpoint{0.000000in}{-0.500000in}}{\pgfqpoint{9.000000in}{0.3500000in}}%
\pgfusepath{use as bounding box}%
\begin{pgfscope}%
\pgfsetbuttcap%
\pgfsetroundjoin%
\definecolor{currentfill}{rgb}{0.300000,0.300000,0.300000}%
\pgfsetfillcolor{currentfill}%
\pgfsetfillopacity{0.500000}%
\pgfsetlinewidth{1.003750pt}%
\definecolor{currentstroke}{rgb}{0.300000,0.300000,0.300000}%
\pgfsetstrokecolor{currentstroke}%
\pgfsetstrokeopacity{0.500000}%
\pgfsetdash{}{0pt}%
\pgfpathmoveto{\pgfqpoint{0.070378in}{-0.812929in}}%
\pgfpathlineto{\pgfqpoint{9.210178in}{-0.812929in}}%
\pgfpathlineto{\pgfqpoint{9.210178in}{-0.058611in}}%
\pgfpathlineto{\pgfqpoint{0.070378in}{-0.058611in}}%
\pgfpathlineto{\pgfqpoint{0.070378in}{-0.812929in}}%
\pgfpathclose%
\pgfusepath{stroke,fill}%
\end{pgfscope}%
\begin{pgfscope}%
\pgfsetbuttcap%
\pgfsetroundjoin%
\definecolor{currentfill}{rgb}{1.000000,1.000000,1.000000}%
\pgfsetfillcolor{currentfill}%
\pgfsetlinewidth{1.003750pt}%
\definecolor{currentstroke}{rgb}{0.000000,0.000000,0.000000}%
\pgfsetstrokecolor{currentstroke}%
\pgfsetdash{}{0pt}%
\pgfpathmoveto{\pgfqpoint{0.042600in}{-0.785151in}}%
\pgfpathlineto{\pgfqpoint{9.182400in}{-0.785151in}}%
\pgfpathlineto{\pgfqpoint{9.182400in}{-0.030833in}}%
\pgfpathlineto{\pgfqpoint{0.042600in}{-0.030833in}}%
\pgfpathlineto{\pgfqpoint{0.042600in}{-0.785151in}}%
\pgfpathclose%
\pgfusepath{stroke,fill}%
\end{pgfscope}%
\begin{pgfscope}%
\pgfsetrectcap%
\pgfsetroundjoin%
\pgfsetlinewidth{1.505625pt}%
\definecolor{currentstroke}{rgb}{0.000000,0.000000,1.000000}%
\pgfsetstrokecolor{currentstroke}%
\pgfsetdash{}{0pt}%
\pgfpathmoveto{\pgfqpoint{0.159266in}{-0.250213in}}%
\pgfpathlineto{\pgfqpoint{0.392600in}{-0.250213in}}%
\pgfusepath{stroke}%
\end{pgfscope}%
\begin{pgfscope}%
\pgfsetbuttcap%
\pgfsetroundjoin%
\definecolor{currentfill}{rgb}{0.000000,0.000000,1.000000}%
\pgfsetfillcolor{currentfill}%
\pgfsetlinewidth{0.501875pt}%
\definecolor{currentstroke}{rgb}{0.000000,0.000000,0.000000}%
\pgfsetstrokecolor{currentstroke}%
\pgfsetdash{}{0pt}%
\pgfsys@defobject{currentmarker}{\pgfqpoint{-0.055556in}{-0.055556in}}{\pgfqpoint{0.055556in}{0.055556in}}{%
\pgfpathmoveto{\pgfqpoint{0.000000in}{-0.055556in}}%
\pgfpathcurveto{\pgfqpoint{0.014734in}{-0.055556in}}{\pgfqpoint{0.028866in}{-0.049702in}}{\pgfqpoint{0.039284in}{-0.039284in}}%
\pgfpathcurveto{\pgfqpoint{0.049702in}{-0.028866in}}{\pgfqpoint{0.055556in}{-0.014734in}}{\pgfqpoint{0.055556in}{0.000000in}}%
\pgfpathcurveto{\pgfqpoint{0.055556in}{0.014734in}}{\pgfqpoint{0.049702in}{0.028866in}}{\pgfqpoint{0.039284in}{0.039284in}}%
\pgfpathcurveto{\pgfqpoint{0.028866in}{0.049702in}}{\pgfqpoint{0.014734in}{0.055556in}}{\pgfqpoint{0.000000in}{0.055556in}}%
\pgfpathcurveto{\pgfqpoint{-0.014734in}{0.055556in}}{\pgfqpoint{-0.028866in}{0.049702in}}{\pgfqpoint{-0.039284in}{0.039284in}}%
\pgfpathcurveto{\pgfqpoint{-0.049702in}{0.028866in}}{\pgfqpoint{-0.055556in}{0.014734in}}{\pgfqpoint{-0.055556in}{0.000000in}}%
\pgfpathcurveto{\pgfqpoint{-0.055556in}{-0.014734in}}{\pgfqpoint{-0.049702in}{-0.028866in}}{\pgfqpoint{-0.039284in}{-0.039284in}}%
\pgfpathcurveto{\pgfqpoint{-0.028866in}{-0.049702in}}{\pgfqpoint{-0.014734in}{-0.055556in}}{\pgfqpoint{0.000000in}{-0.055556in}}%
\pgfpathclose%
\pgfusepath{stroke,fill}%
}%
\begin{pgfscope}%
\pgfsys@transformshift{0.159266in}{-0.250213in}%
\pgfsys@useobject{currentmarker}{}%
\end{pgfscope}%
\begin{pgfscope}%
\pgfsys@transformshift{0.392600in}{-0.250213in}%
\pgfsys@useobject{currentmarker}{}%
\end{pgfscope}%
\end{pgfscope}%
\begin{pgfscope}%
\pgftext[x=0.575933in,y=-0.308546in,left,base]{{\sffamily\fontsize{20.000000}{24.000000}\selectfont $LU$ factorization (\superlu)}}%
\end{pgfscope}%
\begin{pgfscope}%
\pgfsetrectcap%
\pgfsetroundjoin%
\pgfsetlinewidth{1.505625pt}%
\definecolor{currentstroke}{rgb}{0.000000,0.500000,0.000000}%
\pgfsetstrokecolor{currentstroke}%
\pgfsetdash{}{0pt}%
\pgfpathmoveto{\pgfqpoint{0.159266in}{-0.602372in}}%
\pgfpathlineto{\pgfqpoint{0.392600in}{-0.602372in}}%
\pgfusepath{stroke}%
\end{pgfscope}%
\begin{pgfscope}%
\pgfsetbuttcap%
\pgfsetmiterjoin%
\definecolor{currentfill}{rgb}{0.000000,0.500000,0.000000}%
\pgfsetfillcolor{currentfill}%
\pgfsetlinewidth{0.501875pt}%
\definecolor{currentstroke}{rgb}{0.000000,0.000000,0.000000}%
\pgfsetstrokecolor{currentstroke}%
\pgfsetdash{}{0pt}%
\pgfsys@defobject{currentmarker}{\pgfqpoint{-0.055556in}{-0.055556in}}{\pgfqpoint{0.055556in}{0.055556in}}{%
\pgfpathmoveto{\pgfqpoint{-0.055556in}{-0.055556in}}%
\pgfpathlineto{\pgfqpoint{0.055556in}{-0.055556in}}%
\pgfpathlineto{\pgfqpoint{0.055556in}{0.055556in}}%
\pgfpathlineto{\pgfqpoint{-0.055556in}{0.055556in}}%
\pgfpathclose%
\pgfusepath{stroke,fill}%
}%
\begin{pgfscope}%
\pgfsys@transformshift{0.159266in}{-0.602372in}%
\pgfsys@useobject{currentmarker}{}%
\end{pgfscope}%
\begin{pgfscope}%
\pgfsys@transformshift{0.392600in}{-0.602372in}%
\pgfsys@useobject{currentmarker}{}%
\end{pgfscope}%
\end{pgfscope}%
\begin{pgfscope}%
\pgftext[x=0.575933in,y=-0.660705in,left,base]{{\sffamily\fontsize{20.000000}{24.000000}\selectfont Selected Inversion (\pselinv)}}%
\end{pgfscope}%
\begin{pgfscope}%
\pgfsetbuttcap%
\pgfsetroundjoin%
\pgfsetlinewidth{1.505625pt}%
\definecolor{currentstroke}{rgb}{1.000000,0.000000,0.000000}%
\pgfsetstrokecolor{currentstroke}%
\pgfsetdash{{6.000000pt}{6.000000pt}}{0.000000pt}%
\pgfpathmoveto{\pgfqpoint{4.818859in}{-0.250213in}}%
\pgfpathlineto{\pgfqpoint{5.052192in}{-0.250213in}}%
\pgfusepath{stroke}%
\end{pgfscope}%
\begin{pgfscope}%
\pgfsetbuttcap%
\pgfsetroundjoin%
\definecolor{currentfill}{rgb}{1.000000,0.000000,0.000000}%
\pgfsetfillcolor{currentfill}%
\pgfsetlinewidth{0.501875pt}%
\definecolor{currentstroke}{rgb}{0.000000,0.000000,0.000000}%
\pgfsetstrokecolor{currentstroke}%
\pgfsetdash{}{0pt}%
\pgfsys@defobject{currentmarker}{\pgfqpoint{-0.055556in}{-0.055556in}}{\pgfqpoint{0.055556in}{0.055556in}}{%
\pgfpathmoveto{\pgfqpoint{0.000000in}{-0.055556in}}%
\pgfpathcurveto{\pgfqpoint{0.014734in}{-0.055556in}}{\pgfqpoint{0.028866in}{-0.049702in}}{\pgfqpoint{0.039284in}{-0.039284in}}%
\pgfpathcurveto{\pgfqpoint{0.049702in}{-0.028866in}}{\pgfqpoint{0.055556in}{-0.014734in}}{\pgfqpoint{0.055556in}{0.000000in}}%
\pgfpathcurveto{\pgfqpoint{0.055556in}{0.014734in}}{\pgfqpoint{0.049702in}{0.028866in}}{\pgfqpoint{0.039284in}{0.039284in}}%
\pgfpathcurveto{\pgfqpoint{0.028866in}{0.049702in}}{\pgfqpoint{0.014734in}{0.055556in}}{\pgfqpoint{0.000000in}{0.055556in}}%
\pgfpathcurveto{\pgfqpoint{-0.014734in}{0.055556in}}{\pgfqpoint{-0.028866in}{0.049702in}}{\pgfqpoint{-0.039284in}{0.039284in}}%
\pgfpathcurveto{\pgfqpoint{-0.049702in}{0.028866in}}{\pgfqpoint{-0.055556in}{0.014734in}}{\pgfqpoint{-0.055556in}{0.000000in}}%
\pgfpathcurveto{\pgfqpoint{-0.055556in}{-0.014734in}}{\pgfqpoint{-0.049702in}{-0.028866in}}{\pgfqpoint{-0.039284in}{-0.039284in}}%
\pgfpathcurveto{\pgfqpoint{-0.028866in}{-0.049702in}}{\pgfqpoint{-0.014734in}{-0.055556in}}{\pgfqpoint{0.000000in}{-0.055556in}}%
\pgfpathclose%
\pgfusepath{stroke,fill}%
}%
\begin{pgfscope}%
\pgfsys@transformshift{4.818859in}{-0.250213in}%
\pgfsys@useobject{currentmarker}{}%
\end{pgfscope}%
\begin{pgfscope}%
\pgfsys@transformshift{5.052192in}{-0.250213in}%
\pgfsys@useobject{currentmarker}{}%
\end{pgfscope}%
\end{pgfscope}%
\begin{pgfscope}%
\pgftext[x=5.235525in,y=-0.308546in,left,base]{{\sffamily\fontsize{20.000000}{24.000000}\selectfont \REVTWO{$LDL^T$} factorization (\mumps)}}%
\end{pgfscope}%
\begin{pgfscope}%
\pgfsetbuttcap%
\pgfsetroundjoin%
\pgfsetlinewidth{1.505625pt}%
\definecolor{currentstroke}{rgb}{0.000000,0.750000,0.750000}%
\pgfsetstrokecolor{currentstroke}%
\pgfsetdash{{6.000000pt}{6.000000pt}}{0.000000pt}%
\pgfpathmoveto{\pgfqpoint{4.818859in}{-0.602372in}}%
\pgfpathlineto{\pgfqpoint{5.052192in}{-0.602372in}}%
\pgfusepath{stroke}%
\end{pgfscope}%
\begin{pgfscope}%
\pgfsetbuttcap%
\pgfsetmiterjoin%
\definecolor{currentfill}{rgb}{0.000000,0.750000,0.750000}%
\pgfsetfillcolor{currentfill}%
\pgfsetlinewidth{0.501875pt}%
\definecolor{currentstroke}{rgb}{0.000000,0.000000,0.000000}%
\pgfsetstrokecolor{currentstroke}%
\pgfsetdash{}{0pt}%
\pgfsys@defobject{currentmarker}{\pgfqpoint{-0.055556in}{-0.055556in}}{\pgfqpoint{0.055556in}{0.055556in}}{%
\pgfpathmoveto{\pgfqpoint{-0.055556in}{-0.055556in}}%
\pgfpathlineto{\pgfqpoint{0.055556in}{-0.055556in}}%
\pgfpathlineto{\pgfqpoint{0.055556in}{0.055556in}}%
\pgfpathlineto{\pgfqpoint{-0.055556in}{0.055556in}}%
\pgfpathclose%
\pgfusepath{stroke,fill}%
}%
\begin{pgfscope}%
\pgfsys@transformshift{4.818859in}{-0.602372in}%
\pgfsys@useobject{currentmarker}{}%
\end{pgfscope}%
\begin{pgfscope}%
\pgfsys@transformshift{5.052192in}{-0.602372in}%
\pgfsys@useobject{currentmarker}{}%
\end{pgfscope}%
\end{pgfscope}%
\begin{pgfscope}%
\pgftext[x=5.235525in,y=-0.660705in,left,base]{{\sffamily\fontsize{20.000000}{24.000000}\selectfont Selected Inversion (\mumps)}}%
\end{pgfscope}%
\end{tikzpicture}%
\makeatother%
\endgroup%
\end{adjustbox}
}

\subfloat[N=7752, NNZ=2430642]{
\begin{adjustbox}{width=.315\linewidth}
\input{ss_Mumps_LU_DNA.pgf}
\end{adjustbox}
}
\REVTWO{
\subfloat[N=16000, NNZ=12800000]{
\begin{adjustbox}{width=.315\linewidth}
\input{ss_Mumps_DGDFT_ACPNR4_60.pgf}
\end{adjustbox}
}
}
\subfloat[N=32916, NNZ=23857418]{
\begin{adjustbox}{width=.315\linewidth}
\input{ss_Mumps_LU_C_BN_C_1by1.pgf}
\end{adjustbox}
}

\caption{\REV{Performance comparison against \mumps 4.10.0 for computing
selected elements of the inverse.}}
\label{fig:perf_mumps}
\end{figure}

}

Overall, the strong scalability of \pselinv is similar to that of
\superlu \REV{~and it clearly outperforms current inversion algorithm as
implemented in \mumps}.
It requires a modest amount of additional memory to 
compute the selected elements of the inverse. This observation
demonstrates both the validity of our approach and the efficiency of 
our implementation.  More importantly, on matrices arising from actual
electronic calculations, \pselinv allows one to use thousands of \REVTWO{cores},
thereby enabling very large-scale computations.

\subsection{The Accuracy of \pselinv}\label{subsec:accuracy} 

In exact arithmetic, the selected inversion method is an exact
method for computing the selected elements of $A^{-1}$, regardless of
whether $A$ is positive definite or not.  In practice, the selected
inversion method cannot give an exact result due to the presence of
round off errors.  
For sparse direct solver, dynamical pivoting strategies such as the
Bunch-Kaufman process~\cite{BunchKaufman1977} for $LDL^{T}$
factorization has been shown to be effective for reducing the numerical
error especially for indefinite matrices.  On distributed memory
machines, dynamic pivoting strategies can significantly affect the
load balance and the scalability of the factorization process. Thus, they are
not used in \superlu~\cite{LiDemmel2003}.

Since the primary goal of the current implementation of the \pselinv
method is to achieve high parallel scalability, we choose not to 
perform additional pivoting steps after the
factorization step.  We show that for the test problems we tried,
\pselinv \REV{~produces results comparable to that produced by the 
\mumps package, and}
is sufficiently accurate even when the matrix is relatively ill
conditioned.

\REV{The first set of experiments reports comparative results between \mumps and \pselinv
on two problems benchmarked in Section~\ref{subsec:scalability}. 
Selected elements are computed with both methods and their maximum relative
column-wise difference \linebreak $\max_{1\leq j \leq n} ( {\norm{ A^{-1}_{\mumps_{*,j}} -  A^{-1}_{\pselinv_{*,j}} } }/{\norm{ A^{-1}_{\mumps_{*,j}} }} ) $
is presented in Table~\ref{tab:mumps_err}.
In both cases, \pselinv provides accurate results which are comparable to \mumps on these two problems.

\begin{table}
\centering
\REV{
\subfloat[\ludna]{
\begin{tabular}{|c|c|}
\hline
$P$ & Max. relative error \\
\hline
1   &   8.52E-008\\
4   &   2.35E-008\\
16  &   2.52E-008\\
36  &   2.05E-007\\
64  &   2.07E-007\\
\hline
\end{tabular}
}
\subfloat[\lucbnc]{
\begin{tabular}{|c|c|}
\hline
$P$ & Max. relative error\\
\hline
16     & 1.96E-007\\
36     & 1.72E-007\\
64     & 1.77E-007\\
121    & 1.76E-007\\
256    & 1.78E-007\\
\hline
\end{tabular}
}
}
\caption{\REV{Maximum column-wise relative error between \mumps and \pselinv
using different number of \REVTWO{cores} $P$.}}
\label{tab:mumps_err}
\end{table}

}

\REV{In the second set of experiments, we assess the accuracy of \pselinv on
larger ill conditioned matrices. }
To obtain test problems that are indefinite and ill conditioned, 
for each test problem listed in Tab.~\ref{tab:characteristics},
we construct a sequence of $A(z)$ defined by Eq.~\eqref{eq:az} 
for a number of complex shifts $z$.  The  real parts of the shifts lie
within the spectrum of the matrix pencil $(H,S)$, and the imaginary
parts range from small ($10^{-7}$) to large ($10^{-1}$) values. 

In order to quantify the accuracy of the \pselinv method for large matrices 
in the test set, we need to find an appropriate error metric.
\REV{Due to the large matrix size, the parallel matrix inversion method
in \mumps becomes too expensive to be used for the purpose of benchmarking.}
Motivated by the trace estimation in Eq.~\eqref{eqn:traceest}, we choose to measure the numerical 
error introduced by \pselinv by using the following quantity
\begin{equation}
  E(z) = \frac{\abs{N - \Tr[A(z)^{-1} A(z)]}}{N} 
  \equiv \left| 1 - \frac{1}{N}\sum_{i,j=1}^{N} [A(z)^{-1}]_{ij}
  [A(z)]_{ji}\right|.
  \label{eq:trcerr}
\end{equation}
Since $A$ is sparse, the summation in Eq.~\eqref{eq:trcerr} involves
only those $i$ and $j$ such that $A_{i,j}\neq 0$.

\begin{figure}[h]
	\begin{center}
		\subfloat[Spectral density]{\label{fig:pwtkAccuracy.dos}\includegraphics[width=0.48\textwidth]{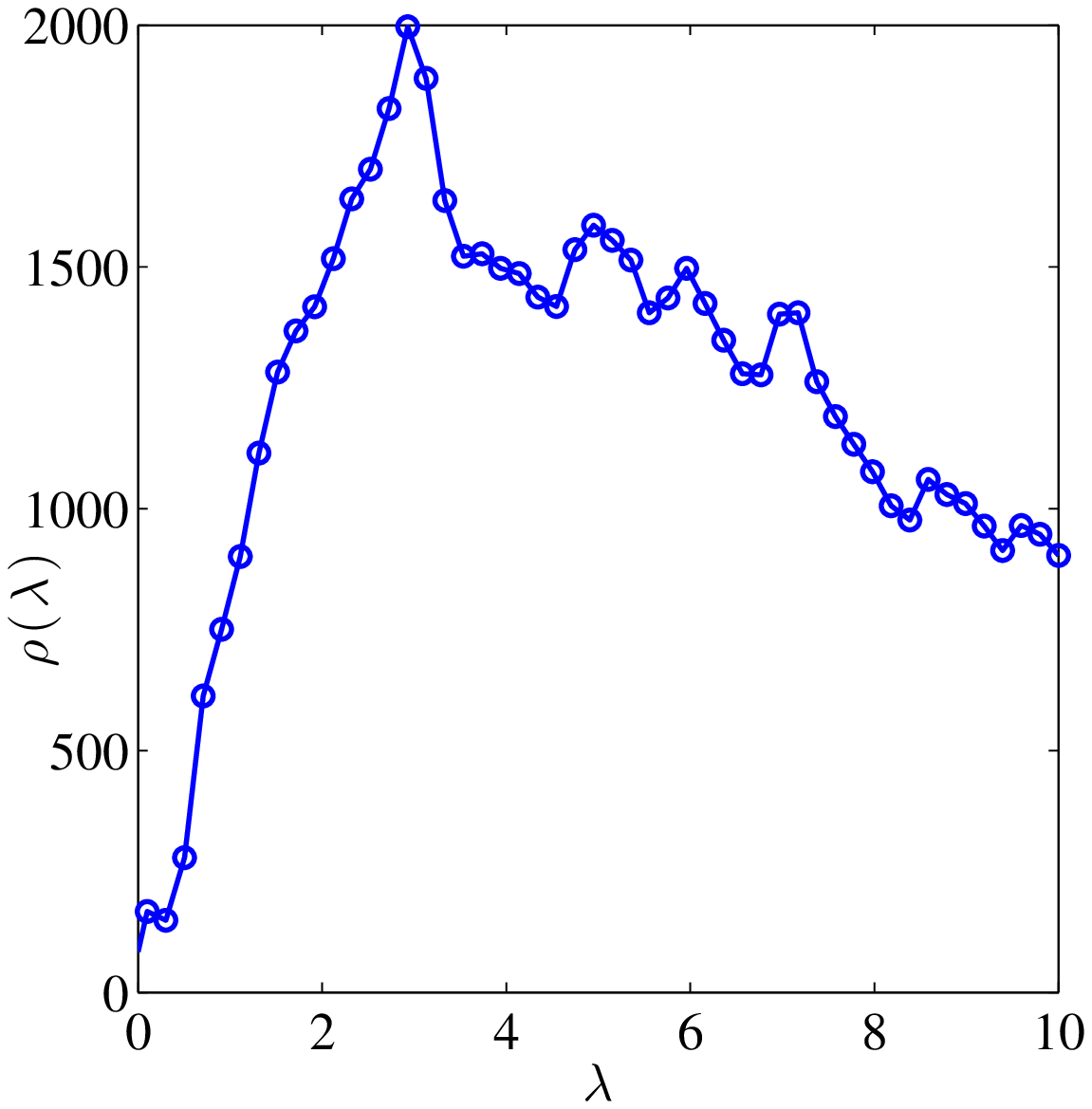}}
		\subfloat[Error $E(z)$ for shifts with different real and imaginary parts.]{\label{fig:pwtkAccuracy.error}\includegraphics[width=0.48\textwidth]{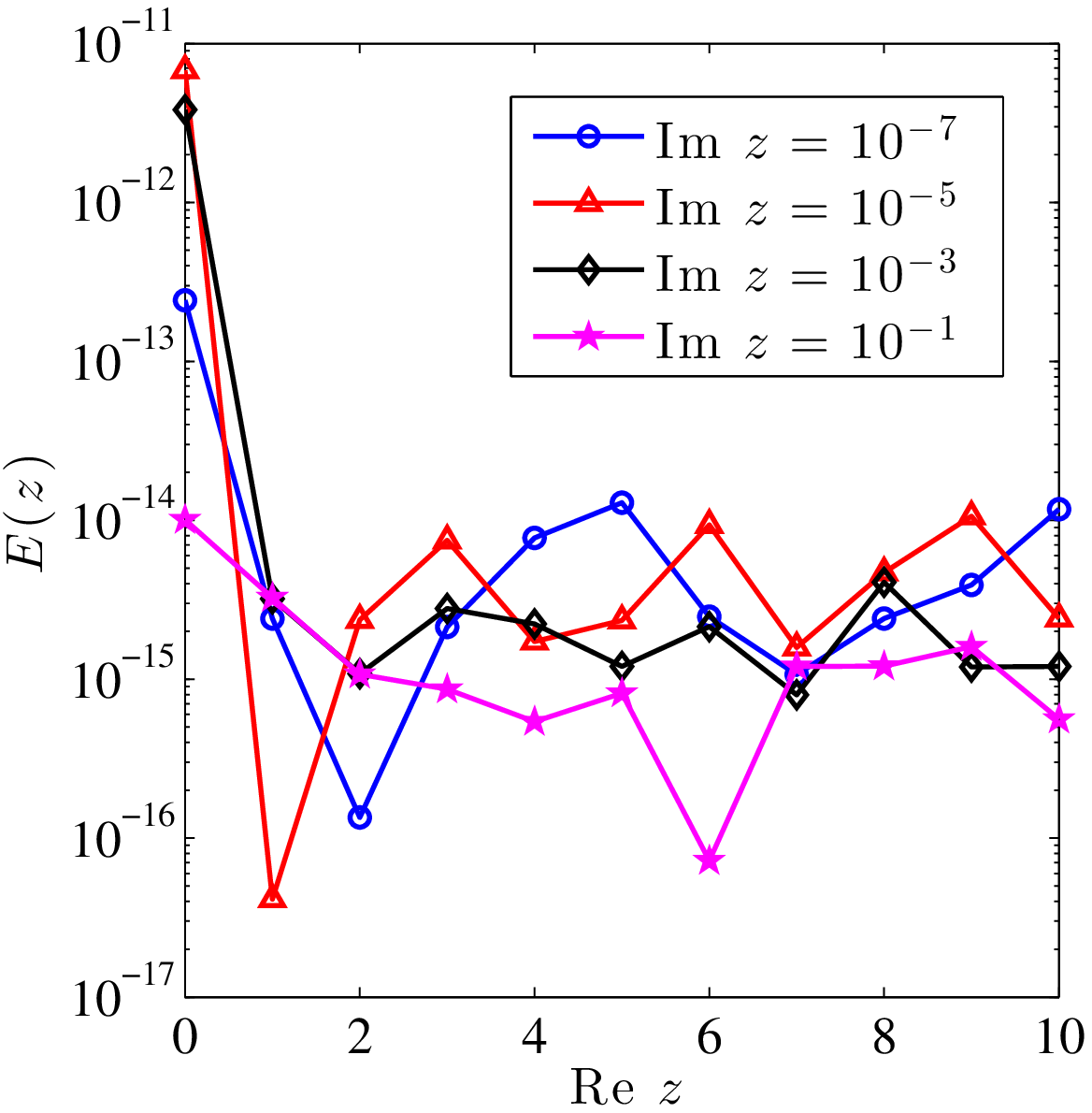}}
	\end{center}
	\caption{Spectral density and error $E(z)$ for the pwtk matrix.}
	\label{fig:pwtkAccuracy}
\end{figure}

\begin{figure}[h]
	\begin{center}
		\subfloat[Spectral density]{\label{fig:DNA64cellAccuracy.dos}\includegraphics[width=0.48\textwidth]{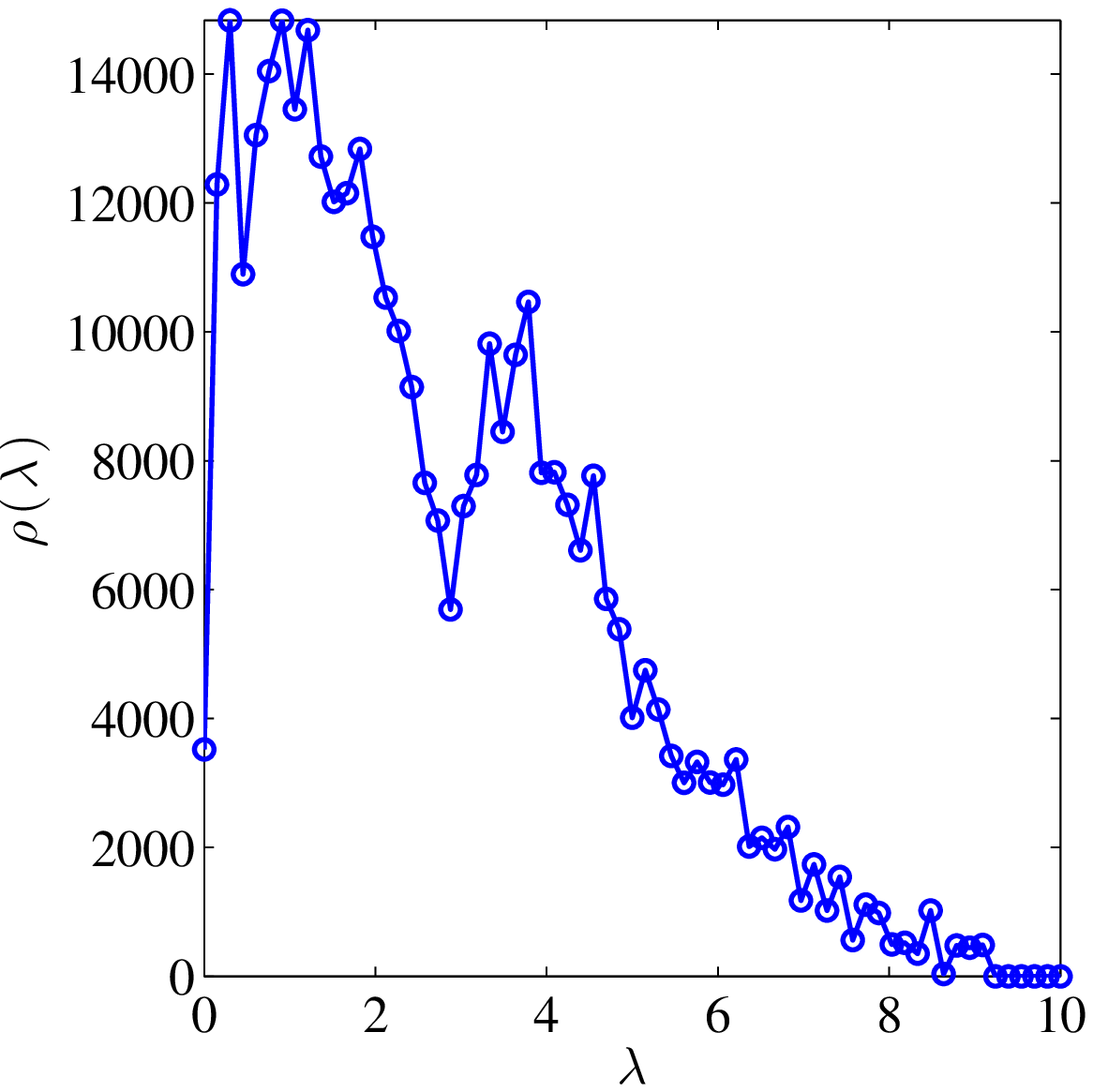}}
		\subfloat[Error $E(z)$ for shifts with different real and imaginary parts.]{\label{fig:DNA64cellAccuracy.error}\includegraphics[width=0.48\textwidth]{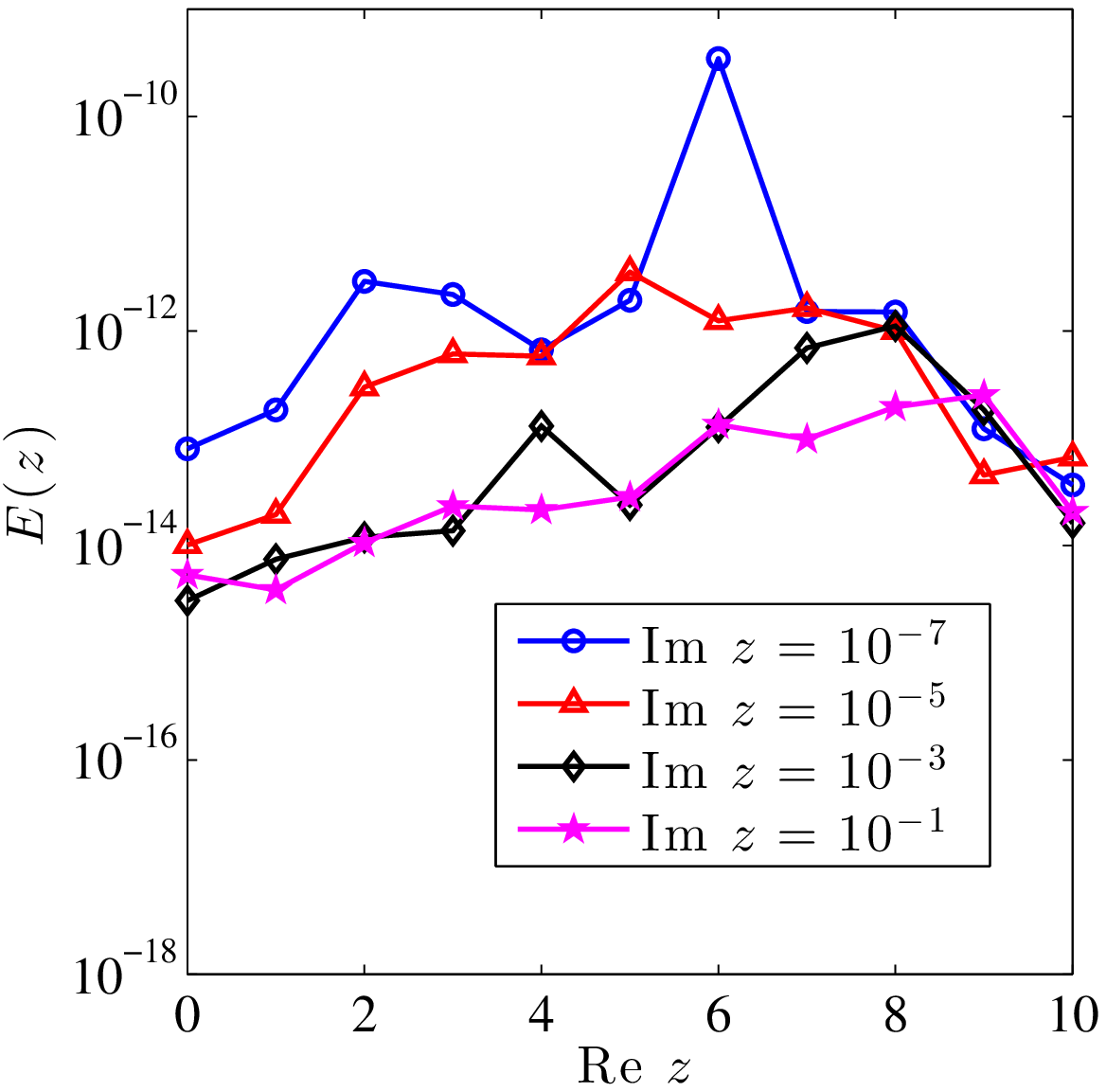}}
	\end{center}
	\caption{Spectral density and error $E(z)$ for the DNA\_715\_64cell matrix.}
	\label{fig:DNA64cellAccuracy}
\end{figure}

In Fig.~\ref{fig:pwtkAccuracy}, we show both the spectral density 
$\rho(\lambda)$ of the pwtk matrix, which describes the number 
of eigenvalues of eigenvalues $(H,S)$ per unit interval, 
and $E(z)$ for a number of shifts $z$ with different real and imaginary
parts.  
We observe that for all these problems, the measured 
errors are below $10^{-11}$, even when the real part of $z$ is close 
to an eigenvalue cluster and the imaginary part of $z$ is as
small as $10^{-7}$.  Fig.~\ref{fig:DNA64cellAccuracy} shows 
that a similar level of accuracy is achieved in the test of
DNA\_715\_64cell matrix.

\subsection{Application to electronic structure theory}\label{subsec:ksdft}

In this section, we demonstrate how \pselinv can be applied to
accelerate Kohn-Sham density functional theory (KSDFT)
calculation~\cite{HohenbergKohn1964,KohnSham1965}, which is widely used
for describing the ground state electronic
properties of molecules, solids and other nano structures. We use the
recently developed pole expansion and selected inversion technique
(PEXSI)~\REV{[Lin et al. 2009a; 2009b; 2011b;
\citeyearNP{LinChenYangEtAl2013}]}
to compute the non-zero elements of the so-called single particle 
density matrix $\Gamma$ that can be approximated by
\begin{equation}
	\Gamma \approx \sum_{l=1}^{P} \omega_l(H-z_{l}S)^{-1}.
	\label{eqn:expandpole}
\end{equation}
where $z_{l},\omega_{l}\in \mathbb{C}$, and the number of ``poles''
$P$ is around $80$ in practical calculations. 
Both $H$ and $S$ are real sparse symmetric matrices that have the same 
sparsity pattern when a local basis set is used to discretize
the Kohn-Sham problem. Each matrix $A_{l}=H-z_{l}S$ is a
complex symmetric matrix. 
\REV{In KSDFT calculations, only elements of $\Gamma_{ij}$ corresponding to
the nonzero elements of $H$ and $S$ (i.e. $H_{ij},S_{ij}\ne 0$) are
needed to compute physical quantities such as electron density and
energy. The expansion~\eqref{eqn:expandpole} immediately suggests that
only the selected elements as defined in Eq.~\eqref{eqn:selelem} of the
matrices $A_{l}^{-1}$ are needed.}
In a parallel implementation of PEXSI, 
we use \pselinv to evaluate the selected elements of $A_{l}^{-1}$ 
that correspond to the nonzero elements of $\Gamma$ on a subset of 
\REVTWO{cores}. The selected inversion of $A_{l}$ for different $l$
can be carried out independently on different subsets of \REVTWO{cores}.

We apply the parallel PEXSI method to the DG\_Graphene\_2048 and
\linebreak DG\_Graphene\_8192 systems, which are disordered graphene 
systems with $2048$ and $8192$ atoms, respectively, and compare its
performance with a standard approach that requires a partial 
diagonalization of $(H,S)$. We use a ScaLAPACK subroutine
\textsf{pdsyevr}~\cite{Vomel2010}, which is based on the multiple 
relatively robust representations (MRRR) algorithm, to perform
such diagonalization. \REV{Though both $H$ and $S$ are sparse
matrices, the MRRR algorithm treats them as dense matrices. For $H,S\in
\mathbb{R}^{N\times N}$, the MRRR algorithm first performs a
tridiagonalization procedure with $\Or(N^3)$ cost, then efficiently
solves the eigenvalues and eigenvectors of the tridiagonal system with
$\Or(N^2)$ cost, and finally the eigenvectors can be constructed with
$\Or(N^3)$ cost.} 

\begin{figure}[h]
	\begin{center}
		\subfloat[$2048$ atoms graphene system.]{\label{fig:timeGraphene2048_8192.2048}\includegraphics[width=0.48\textwidth]{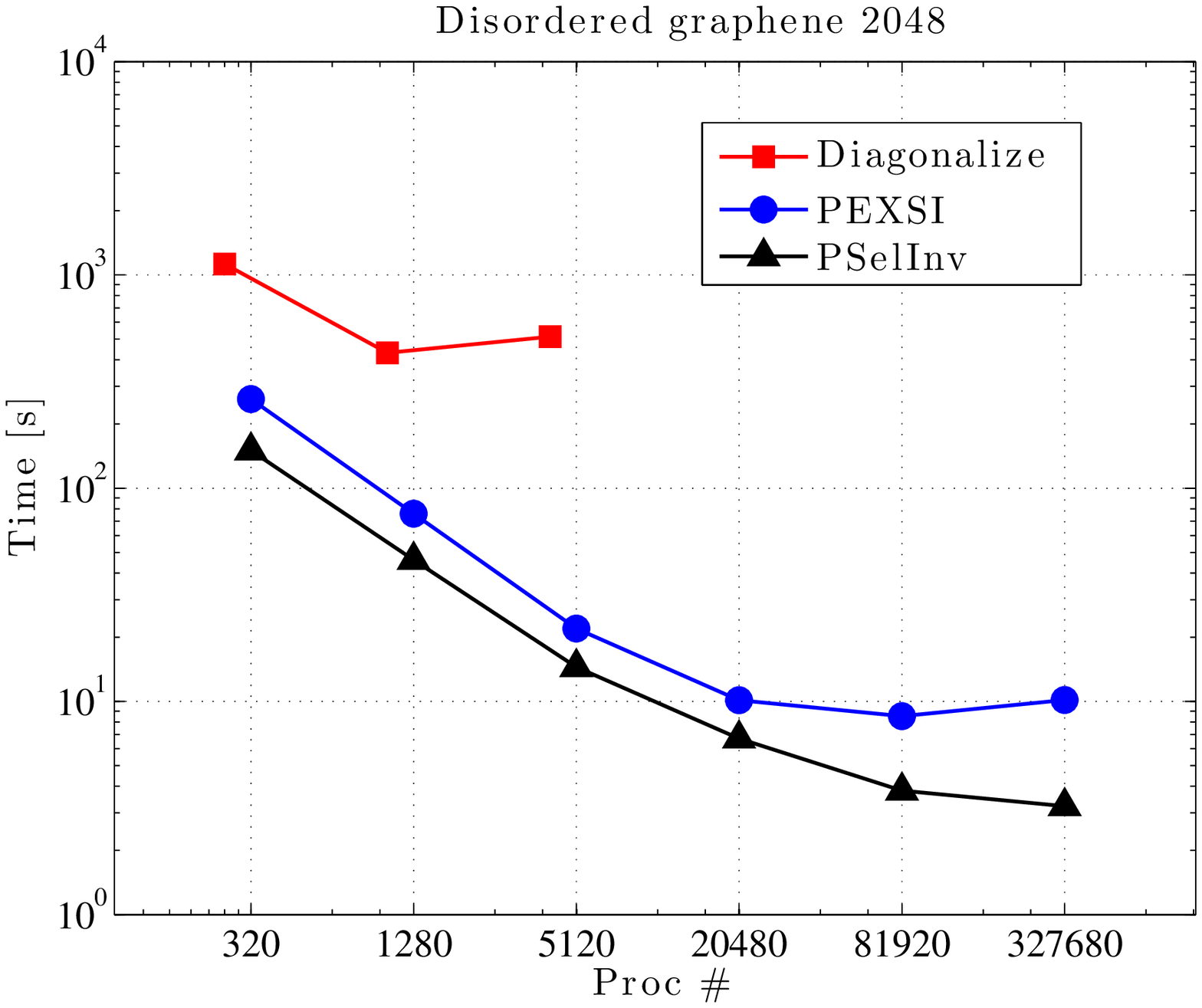}}
		\subfloat[$8192$ atoms graphene system.]{\label{fig:timeGraphene2048_8192.8192}\includegraphics[width=0.48\textwidth]{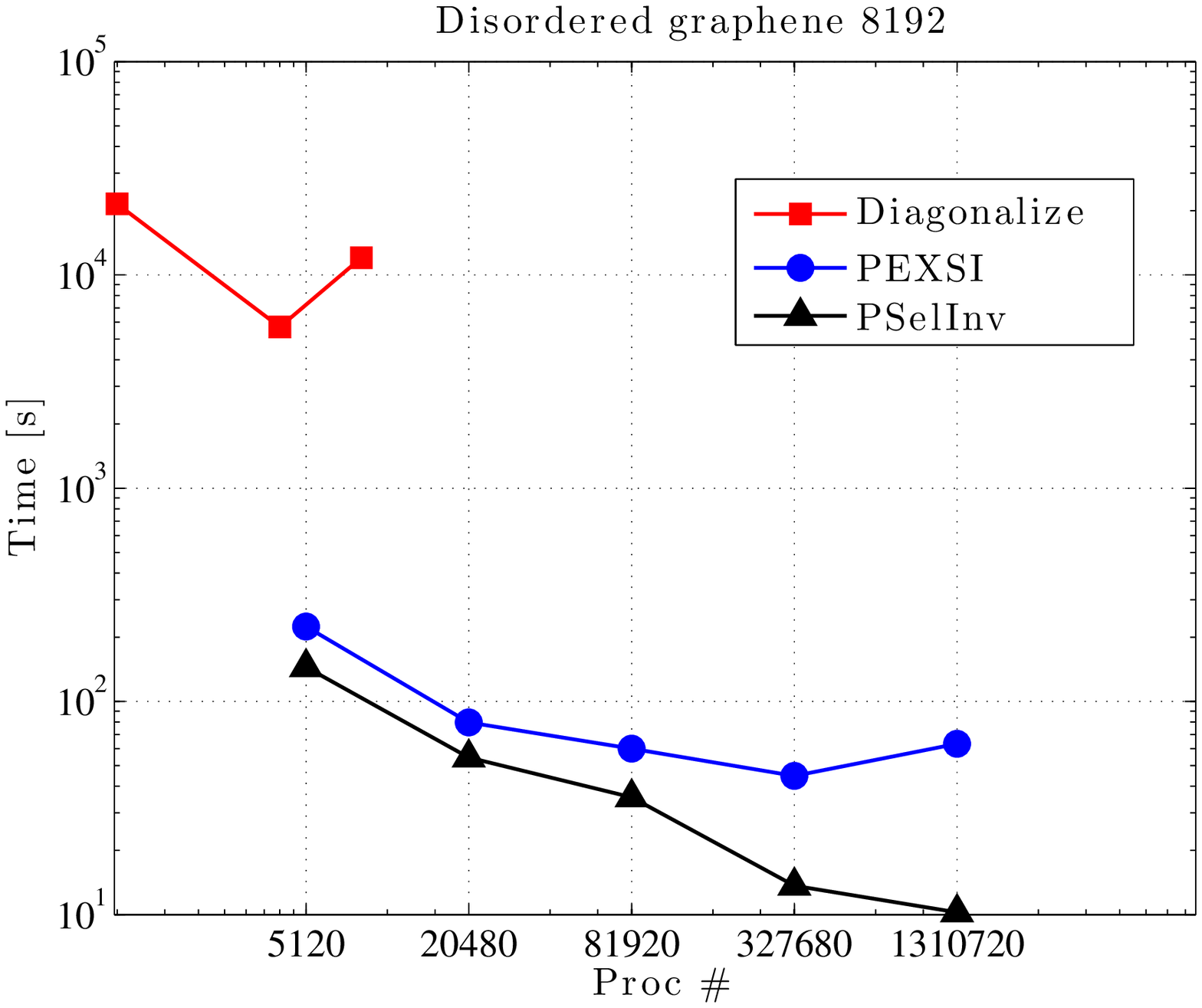}}
	\end{center}
	\caption{The wall clock time versus the number of \REVTWO{cores} for a graphene system.  }
	\label{fig:timeGraphene2048_8192}
\end{figure}

Fig.~\subref*{fig:timeGraphene2048_8192.2048} shows that
the scalability of \textsf{pdsyevr} is limited to $1024$ \REVTWO{cores} 
for the $2048$-atom problem.  Adding more \REVTWO{cores} to the diagonalization 
process leads to an increase of the wall clock time due to communication 
overhead. For this relatively small problem, the benefit of parallel
implementation of PEXSI is already clear when 320 \REVTWO{cores} are used.
Since we use $P=80$ poles in the pole expansion, $4$ \REVTWO{cores}
\REV{with a $2\times 2$ processor grid} are used
in each selected inversion in this case. The wall clock time used by PEXSI 
is $261$ seconds, among which $150$ seconds are attributed to \pselinv,
\REV{$95$ seconds are attributed to factorization, and $10$ seconds for
symbolic factorization}.  
This timing result compares favorably to the $430$ seconds of 
measured wall clock time required by \textsf{pdsyevr} on 1,024 \REVTWO{cores}.

Furthermore, we can clearly see from Fig.~\subref*{fig:timeGraphene2048_8192.2048} 
that the parallel PEXSI method can scale to a much larger number of \REVTWO{cores}.
Nearly perfect speedup can be observed when a total of $20,480$ \REVTWO{cores}
are used in the parallel PEXSI computation, with $256$ \REVTWO{cores} used in
each selected inversion. The total wall clock time used in this calculation
is merely 10 seconds.  Compared to the best wall clock time we can obtain
for the diagonalization procedure, which is $430$ sec on 1,024 \REVTWO{cores}, 
this represents a speedup factor of $43$.

For the larger system that contains $8192$ atoms, \textsf{pdsyevr} can
scale to 4,096 \REVTWO{cores} as we can see in
Fig.~\subref*{fig:timeGraphene2048_8192.8192}. It takes $5703$ wall
clock seconds to perform such a computation. When the parallel PEXSI
calculation is carried out on $5,120$ \REV{\REVTWO{cores}} with $64$
\REVTWO{cores} used to perform each selected inversion, the total wall clock
time required is $224$ seconds.
Fig.~\subref*{fig:timeGraphene2048_8192.8192} also shows that parallel
PEXSI can scale to as many as $327,680$ \REVTWO{cores} with $4,096$
\REVTWO{cores} used for each selected inversion.  The total wallclock time
required in this calculation is merely $45$ seconds, a 127 fold speedup
compared to the best diagonalization wallclock time measured.

\begin{figure}[h] \begin{center}
	\includegraphics[width=0.65\textwidth]{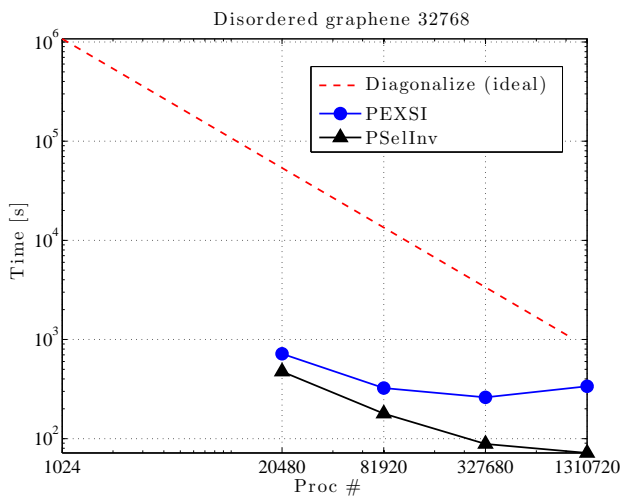}
	\end{center}
	\caption{The wall clock time versus the number of \REVTWO{cores} for
	a graphene system with $32768$ atoms. }
	\label{fig:timeGraphene32768}
\end{figure}

Finally we apply PEXSI to a system with $32,768$ atoms. The matrix size
is $4$ times larger than that for the system with $8,192$ atoms, and the
diagonalization routine is no longer feasible: the wall
clock time required to run \textsf{pdsyevr} routine with $1024$ \REVTWO{cores} for
the $2048$-atom system is $431$ sec, and for the $8192$-atom system is
$21556$ sec. The increase of the wall clock time is $50$ fold, which is
roughly in agreement with the cubic complexity scaling factor $4^3=64$.  
The cubic scaling of the diagonalization procedure implies that
the wall clock time would increase by at least a factor of $50$ to 
$1,077,800$ seconds ($300$ hours) if we perform the same type of calculation 
for a system that contains  $32,768$ atoms on $1,024$ \REVTWO{cores}. Based
on this estimation and assuming that the strong scaling of
\textsf{pdsyevr} is perfect, we compare the ideal performance of the
diagonalization method with the practical performance of PEXSI 
in Fig.~\ref{fig:timeGraphene32768} up to $1,310,720$ \REVTWO{cores}.  The total 
wall clock time for both factorization and
selected inversion reaches its minimum at $4,096$ \REVTWO{cores} per pole
($327,680$ \REVTWO{cores} in total), which is $241$ sec. \REV{Among these,
$87$ sec is attributed to \pselinv.}  Comparatively, even if the 
diagonalization procedure scales perfectly to more than 1 million \REVTWO{cores}, 
which is highly unlikely
within the current framework of diagonalization methods, the projected
wall clock time is over $1000$ sec, which is significantly more than
that used by PEXSI.

\section{Conclusion and future work}\label{sec:conclusion}

We described an efficient parallel implementation of the selected
inversion algorithm for distributed memory parallel machines.  The current
implementation of \pselinv can be applied to the computation of selected
elements of the inverse of a sparse symmetric matrix.  It is publicly
available, and can scale to more than 4000 \REVTWO{cores} for sufficiently
large problems.  The scalability of the solver depends on the size and
sparsity of the matrix. We observed that it is important to exploit
concurrency available within the elimination trees to achieve high
scalability in the parallel selected inversion process.   In the future,
we plan to further improve the tree level parallelism to enhance the
concurrency among different supernodes.  We also observed that, for our
test problems, the \pselinv method is relatively accurate even for matrices
that are highly indefinite and close to singular. It can be applied
to accelerate several scientific computation applications such as the
density functional theory based electronic structure calculations.  
In order to further improve the numerical accuracy of the \pselinv 
method, especially for indefinite matrices, dynamic pivoting strategies 
such as Bunch-Kaufman procedure~\cite{BunchKaufman1977,grls:94} 
for the factorization 
may be needed.  Generalizing \pselinv to non-symmetric matrices and 
combining it with other sparse direct solvers are also areas we plan to
work on in the future.

\begin{acks}
This work was partially supported by the Laboratory Directed Research
and Development Program of Lawrence Berkeley National Laboratory under
the U.S. Department of Energy contract number DE-AC02-05CH11231 (L. L.
and C. Y.), the Scientific Discovery through Advanced Computing (SciDAC)
program funded by U.S. Department of Energy, Office of Science, Advanced
Scientific Computing Research and Basic Energy Sciences (M. J., L. L.
and C. Y.), and the Center for Applied Mathematics for Energy Research
Applications (CAMERA), which is a partnership between Basic Energy
Sciences (BES) and Advanced Scientific Computing Research (ASRC) at the
U.S Department of Energy.  We would like to thank Xiaoye S. Li and
Fran\c{c}ois-Henry Rouet for helpful discussion. 
\end{acks}


\end{document}